\documentclass[12pt]{article}
\usepackage{graphicx}

\newtheorem{theorem}{Theorem}[section]
\newtheorem{corollary}[theorem]{Corollary}

\newtheorem{proposition}[theorem]{Proposition}
\newenvironment{proof}[1][Proof]{\textbf{#1.} }{\
\rule{0.5em}{0.5em}}

\def\text{\hbox} 

\def\a{\alpha}
\def\b{\beta}

\def\l{\lambda}
\def\m{\nu}

\def\p{\pi}

\def\r{\rho}

\def\f{\phi}

\def\o{\omega}
\def\ot{\widetilde{\o}}

\def\S{\Sigma}

\def\T{{\mathcal T}}
\def\K{{\mathcal K}}

\def\Z{{\mathbf Z}}
\def\R{{\mathbf R}}
\def\S{{\mathbf S}}
\def\Di{{\mathbf D}}
\def\bt{{\mathbf b}}
\def\Gt{\widetilde{G}}
\def\St{\widetilde{S}}

\def\T1{T_1({\bf p/q};{\bf r/s})}

\def\K2{{\mathcal K}_2({\bf p/q};{\bf r/s})}

\begin{document}

\title{The many faces of cyclic branched coverings
\\
of 2-bridge knots and links
\thanks{Work performed under the auspices of G.N.S.A.G.A. --
I.N.d.A.M.,
Italy, and supported by the University of Bologna, funds for selected
research topics. The second named author was supported by the Russian
Foundation for Basic Research.}}

\author{Michele Mulazzani \and Andrei Vesnin}

\maketitle

\begin{abstract}
We discuss 3-manifolds which are cyclic coverings of the 3-sphere,
branched over 2-bridge knots and links. Different descriptions of
these manifolds are presented: polyhedral, Heegaard diagram, Dehn
surgery and coloured graph constructions. Using these
descriptions, we give presentations for their fundamental groups,
which are cyclic presentations in the case of 2-bridge knots. The
homology groups are given for a wide class of cases. Moreover, we
prove that each singly-cyclic branched covering of
a 2-bridge link is the composition of a meridian-cyclic branched
covering of a determined link and a cyclic branched
covering of a trivial knot.
\end{abstract}

\bigskip

\centerline{CONTENTS}

1. Introduction and preliminaries.

2. Cyclic branched coverings.

3. Polyhedral construction.

4. Heegaard diagram construction and genus.

5. Surgery construction.

6. Coloured graph construction.

7. Fundamental groups.

8. Homology groups.

9. Decomposition of singly-cyclic branched coverings.

\section{Introduction and preliminaries}

The family of 2-bridge knots/links is a well-studied subject,
starting from the classical papers of Listing, Dehn, Alexander,
Reidemeister, Schubert and others. Almost everything is known
about the symmetries, invariants, and geometric properties of
these links, of their complements, and of manifolds connected with
them. The aim of the present paper is to provide a review of the
properties of 3-manifolds which are cyclic branched coverings of
2-bridge knots/links. We recall that there are a lot of ways of
describing a 3-manifold: by fundamental polyhedra \cite{ST}, by
Heegaard diagrams \cite{He}, by surgery on links in $\S^3$
\cite{Ro}, by branched coverings of $\S^3$ \cite{Al} and by
gems/crystallizations \cite{Li,Pe}. So, each manifold has many
``faces'', which depend on the type of the description we choose.
This choice depends on the nature of the problem one studies. In
this paper we are going to present these aspects of cyclic
branched coverings of 2-bridge knots/links. As an example, we will
describe the different faces of the classical Hantzsche--Wendt
manifold \cite{Zi1}.

By the term {\it manifold} we always mean a compact, connected,
orientable \hbox{PL-manifold} without boundary. Let $M,N$ be
triangulated $m$-dimensional manifolds and let $L$ be an
$(m-2)$-subcomplex of $N$. A non-degenerate map $f: M \to N$ is an
{\it $n$-fold covering map branched over $L$}, with $n>1$, if: (i)
$f_{|M-f^{-1}(L)} : M-f^{-1}(L)\to N-L$ is an ordinary $n$-fold
covering, and (ii) $L=\{x\in N \mid \# f^{-1}(x) < n \}$. The
manifold $M$ is said to be a {\it branched covering} of $N$, and
$L$ is called the {\it branching set} of the covering. Two
branched coverings $f':M'\to N'$ and $f'':M''\to N''$ are
equivalent if there exist two homeomorphisms $\psi : M' \to M''$
and $\phi : N' \to N''$ such that $\psi \circ f'' = f' \circ
\phi$.

A remarkable result by R.~H.~Fox \cite{Fo1} states that a branched
covering is uniquely determined by the ordinary covering induced
by restriction. This proves the existence of a one-to-one
correspondence between the $n$-fold coverings of $N$ branched over
$L$ and the equivalence classes of monodromies (i.e. transitive
representations) $\omega :\pi_1(N-L,x_0)\to\Sigma_n$, where
$\Sigma_n$ is the symmetric group on $n$ elements and $x_0\in N-L$
is an arbitrary base point. If $N$ is an $m$-sphere, then a
covering of $N$ branched over $L$ is simply called a {\it branched
covering of $L$}.

Two $n$-fold branched coverings $f':M'\to N'$, $f'':M''\to N''$,
with branching sets $L'\subset N'$, $L''\subset N''$ and monodromy
maps $\omega_{f'}$, $\omega_{f''}$ respectively, are equivalent if
and only if there exists an inner automorphism $\l$ of $\Sigma_n$
and a homeomorphism $\phi : N' \to N''$, such that $\phi (L') =
L''$ and $\l\circ \omega_{f'} = \omega_{f''}\circ \phi_*$, where
$\phi_*:\p_1(N'-L', x_0')\to\p_1(N''-L'',\phi(x_0'))$ is the
homomorphism induced by $\phi _{|N'-L'}$ on the fundamental
groups.

Branched coverings of spheres are of great interest, in
particular as a method for representing manifolds.  A classical
result concerning this point was obtained by J.~Alexander \cite{Al}
and
states that every $m$-manifold is a covering of $\S^m$, branched
over the $(m-2)$-skeleton of a standard $m$-simplex.

Fox's result gives the possibility of extending the concept of
cyclic coverings from ordinary coverings to branched coverings.
Thus, a branched covering is called {\it cyclic\/} if its
associated ordinary covering is cyclic. Similar extension
also applies to regular and abelian coverings.

Since a cyclic covering is abelian, it is determined up to
equivalence by an epimorphism $$\ot_f:H_1(N-L)\to\Z_n,$$ where
$\Z_n$ is the cyclic group of order $n$, embedded into $\Sigma_n$
through the monomorphism sending $1\in\Z_n$ to the standard cyclic
permutation $(1\,2\,\cdots\,n)\in\Sigma_n$. If $N=\S^3$ and
$L=\bigcup_{j=1}^{\m}L_j$ is a $\m$-component link in $\S^3$, then
$H_1(N-L)\cong\Z^{\m}$ and a basis is given by any set of homology
classes of meridian loops around the components of $L$. Therefore,
an $n$-fold cyclic branched covering $f$ of $L$ is defined by
orienting $L$ and assigning an integer $k_j\in\Z_n-\{0\}$ to each
component $L_j$, such that the set $\{k_1,\ldots,k_{\m}\}$
generates $\Z_n$. If $m_j$ is a meridian around $L_j$, coherently
oriented with the chosen orientations of $L$ and $\S^3$, we define
$\ot_f[m_j]=k_j\in\Z_n$ and therefore
$\o_f[m_j]=(1\,2\,\cdots\,n)^{k_j}$. We will denote this manifold
by $M_{n,k_1,\ldots,k_{\m}}(L)$. By multiplying each $k_j$ by the
same invertible element $u$ of $\Z_n$, we obtain an equivalent
covering. More precisely, two $n$-fold cyclic branched coverings
$f':M'\to N'$ and $f'':M''\to N''$, with associated maps
$\ot_{f'}: H_1(N'-L')\to\Z_n$ and $\ot_{f''}: H_1(N''-L'')\to\Z_n$
respectively, are equivalent if and only if there exist $u\in\Z_n$
and a homeomorphism $\f:N'\to N''$, such that $\gcd(u,n)=1$,
$\f(L')=L''$ and $\ot_{f''}\circ\f_\#=u\ot_{f'}$, where
$\f_\#:H_1(N'-L')\to H_1(N''-L'')$ is the homomorphism induced by
$\f_{\vert N'-L'}$ on the first homology groups and $u\ot_{f'}$ is
the multiplication of $\ot_{f'}$ by $u$. General references on
cyclic branched coverings of knots/links are the interesting books
 \cite{BZ}, \cite{Ka}, and \cite{Ro}.

Following \cite{MM} we shall say that a cyclic branched covering
$M_{n,k_1,\ldots,k_{\m}}(L)$ is:

\begin{description}
\item[(i)] {\it strictly-cyclic\/} if $k_{j'}=k_{j''}$, for every
$j',j''\in\{1,\ldots,\m\}$;
\item[(ii)] {\it almost-strictly-cyclic\/} if $k_{j'}=\pm k_{j''}$,
for
every $j',j''\in\{1,\ldots,\m\}$;
\item[(iii)] {\it meridian-cyclic\/} if $\gcd(n,k_j)=1$, for every
$j\in\{1,\ldots,\m\}$;
\item[(iv)] {\it singly-cyclic\/} if there exists
$j\in\{1,\ldots,\m\}$
such that $\gcd(n,k_j)=1$;
\item[(v)] {\it monodromy-cyclic\/} if it is cyclic.
\end{description}

The following implications are straightforward: $ \text{(i)}
\Rightarrow \text{(ii)} \Rightarrow\text{(iii)} \Rightarrow
\text{(iv)} \Rightarrow \text{(v)}. $ Moreover, all five
definitions are equivalent when either $L$ is a knot or $n=2$.
Note that strictly-cyclic coverings are also called {\it uniform}
in \cite{Zi3} and meridian-cyclic coverings are also called {\it
strongly cyclic} in \cite{Zi4}. By a suitable reorientation of the
link, an almost-strictly-cyclic covering becomes a strictly-cyclic
one. For a singly-cyclic covering we can always assume $k_1=1$, up
to equivalence and possible reordering of the components of $L$.
Therefore, when $\m=2$ the covering is completely determined by an
integer $k=k_2\in\Z_n-\{0\}$.

In order to simplify the notations, the $n$-fold strictly-cyclic
branched covering $M_{n,1,\ldots,1}(L)$ will be denoted by
$M_n(L)$ and the $n$-fold singly-cyclic branched covering of a
2-component link $M_{n,1,k}(L)$ will be denoted by $M_{n,k}(L)$.
In particular, $M_{n,k}(L)$ is meridian-cyclic when $\gcd(n,k)=1$
and strictly-cyclic when $k=1$.

The class of 2-bridge knots/links (also called ``rational
knots/links'') is of great interest and the main properties of its
elements can be found in \cite{BZ, Ka, Ro}. We denote by $\bt
(\a,\b)$ the 2-bridge knot or link of type $(\a,\b)$, with
integers $\a>1$ and $\b\in\Z_{2\a}$ such that $\gcd(\a,\b)=1$. It
is well known that ${\bf b}(\a,\b)$ is a knot when $\a$ is odd and
a two-component link when $\a$ is even. In our discussion we will
skip the ``singular'' case of the 2-component trivial link, which
is a 2-bridge link and is often indicated by $\bt (0,1)$. A
standard diagram presentation for 2-bridge knots/links is
Schubert's normal form \cite[Ch.~12]{BZ}. Actually, Schubert's
normal form is only defined when $\b$ is odd, but this is not a
restriction, up to equivalence (see classification of 2-bridge
knots/links below). Figure~1 gives two examples of such diagrams
for the figure-eight knot $\bt (5,3)$ and for the Whitehead link
$\bt (8,3)$.

\begin{figure}
 \begin{center}
 \includegraphics*[totalheight=16cm]{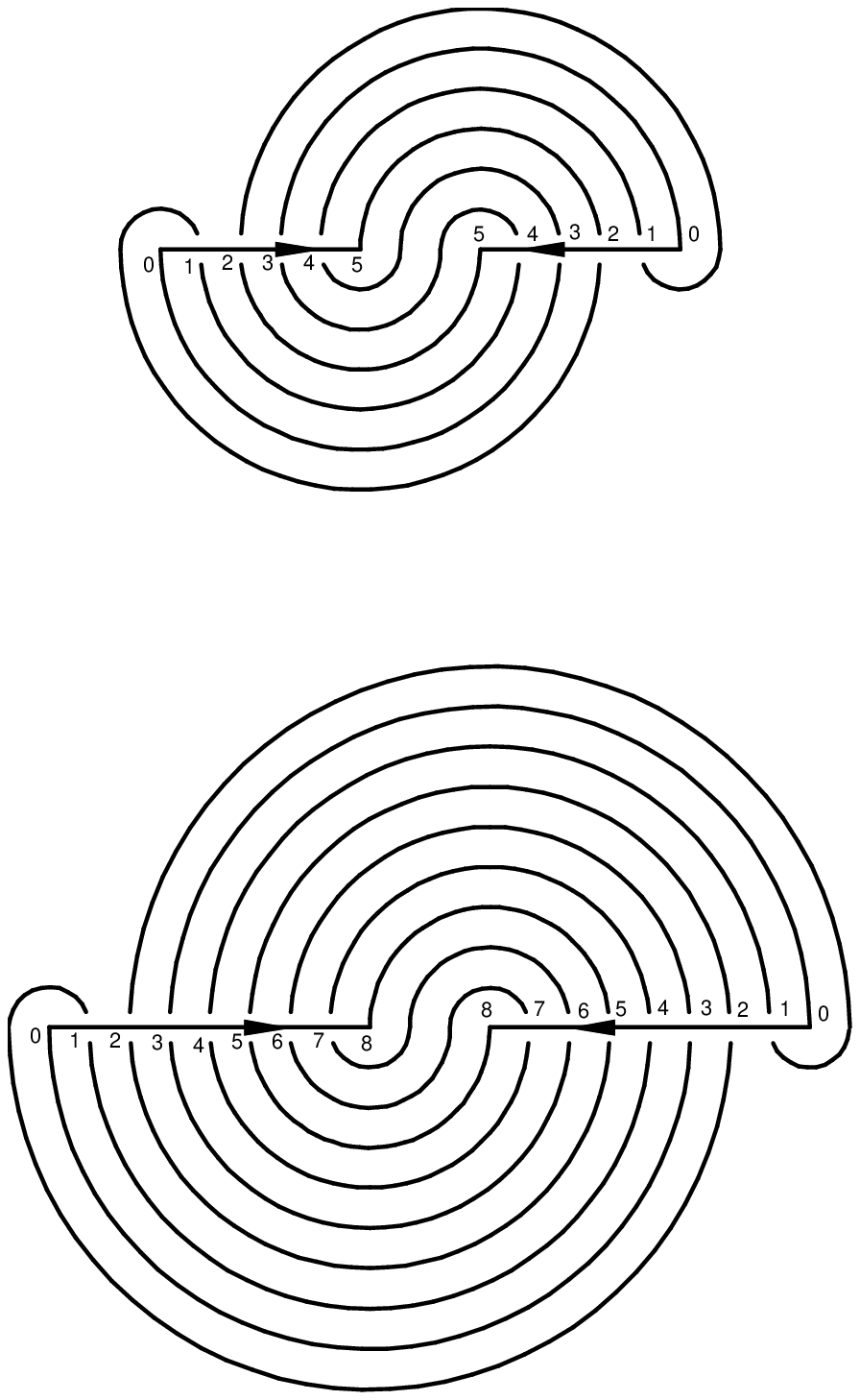}
 \end{center}
 \caption{Schubert's normal form for $\bt (5,3)$ and $\bt (8,3)$.}
 \label{Fig. 1}
\bigskip\bigskip
\end{figure}




Another type of diagram for 2-bridge knots/links is introduced by
J.~Conway \cite{Co} (see Figure~2). It comes from the
representation of a rational number by a continued fraction. If
$$\frac {\a}{\b}  = c_1+ \frac{1}{c_2 +\frac{1} {c_3 + \cdots +
\frac{1}{c_m}}}\,,$$  we will say that $[c_1, c_2, \ldots, c_m]$
are Conway parameters for $\bt (\a,\b)$ and the notation $\a/\b =
[c_1, c_2, \ldots, c_m]$ will be used.

\begin{figure}
 \begin{center}
 \includegraphics*[totalheight=18cm]{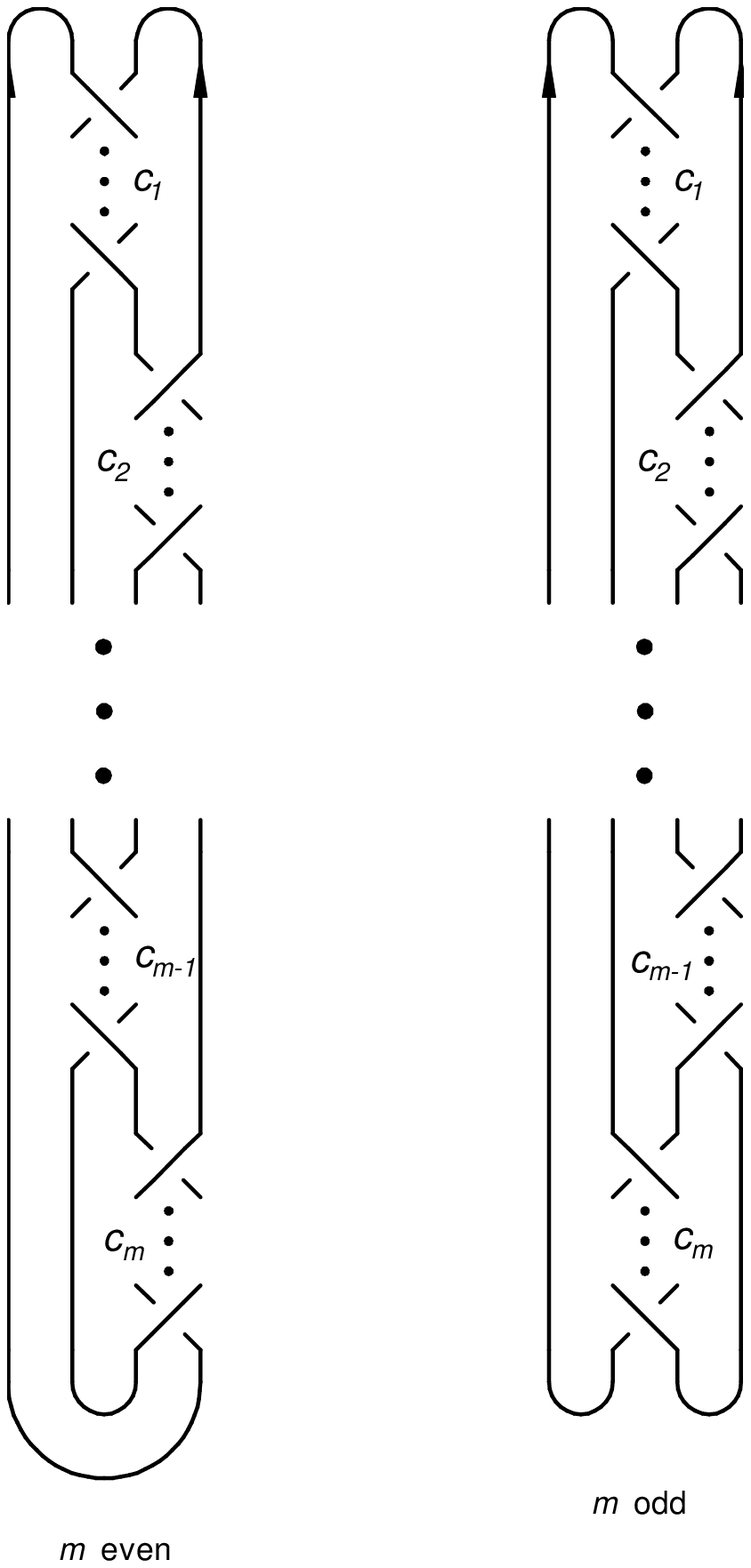}
 \end{center}
 \caption{Conway's normal form for 2-bridge knots and links.}
 \label{Fig. 2}
\bigskip\bigskip
\end{figure}




2-bridge knots/links are assumed to be oriented in a standard way,
as indicated in Figures~1 and 2. By the Schubert construction of
$\bt(\a,\b)$, there is an obvious orientation preserving
involution $\r:\S^3\to\S^3$ exchanging the two bridges of
$\bt(\a,\b)$ and preserving the orientation of the component(s).

Recall that two oriented links $L$ and $L'$ with $\m$ components
with a given order $(L_i)_{1\le i\le \m}$ and $(L'_i)_{1\le i\le
\m}$ are said to be {\it equivalent} if there exists an
orientation preserving homeomorphism $\f : \S^3 \to \S^3$ sending
$L_i$ to $L'_i$, for $i=1,2,\ldots,\m$, and preserving their
orientations. For 2-bridge links the equivalence does not depend
on the ordering of the components because of the homeomorphism
$\r$ previously described. H.~Schubert in \cite{Sc} classified
2-bridge knots/links (see also \cite[Theorem 2.1.3]{Ka}):
\begin{description}
\item [(i)] two oriented two-bridge knots $\bt(\a,\b)$,
$\bt(\a',\b')$
are equivalent if and only if $\a'=\a$ and $\b'\equiv \b^{\pm 1}$
mod $\a$;
\item [(ii)] two oriented two-bridge links $\bt(\a,\b)$,
$\bt(\a',\b')$
are equivalent if and only if $\a'=\a$ and $\b'=\b^{\pm 1}$.
\footnote{Here
and in the following the second coordinate of a 2-bridge knot or link
$\bt(\a,\b)$ will be considered mod $2\a$, except when otherwise
specified.}
\end{description}
If we consider unoriented links, then the condition of (ii)
is reduced to that of (i).

Below we will use the following properties: $\bt(\a,-\b)$ is
equivalent to the mirror image of $\bt(\a,\b)$ and the link
$\bt(\a,\b-\a)$ is equivalent to the link $\bt(\a,\b)$ with the
opposite orientation on one of the two components \cite{BZ}.

{}From now on, we use the notation $M_{n,k',k''}(\a/\b) =
M_{n,k',k''}(\bt(\a,\b))$ for $n$-fold cyclic branched coverings
of 2-bridge links (i.e., $\a$ is even). In
particular, $M_{n,k}(\a/\b)$ will denote a singly-cyclic branched
covering and $M_n(\a/\b)$ will denote a strictly-cyclic one. The
last notation will always be used even in the case when $\bt
(\a,\b)$ is a knot (i.e., $\a$ is odd).

The class of cyclic branched coverings of 2-bridge knots/links has
been intensively studied by many authors, and lot of their
properties were discovered. In particular, the two-fold coverings
are homeomorphic to lens spaces.

In Section~2 we introduce our class of manifolds,
discuss homeomorphisms, geometric structures, and present them as
2-fold branched coverings of $\S^3$.

In Section~3 we will describe the polyhedral construction of these
manifolds, according to J.~Minkus \cite{Mi}.

In Section~4 we will give symmetric Heegaard diagrams for cyclic
branched coverings of 2-bridge knots. Following from a result
obtained in \cite{GM}, this description arises from a construction
by M.~Dunwoody \cite{Du}. In the same section, some estimates for
the genus are given.

In Section~5 we will give a surgery description of these
manifolds, for the case of knots. It was pointed out in \cite{KV}
and \cite{RS} that some cyclic branched coverings of 2-bridge
knots are Takahashi manifolds \cite{Ta} and therefore they can be
obtained by Dehn surgery on a certain chain of circles in $\S^3$.
The Takahashi construction has been generalized in \cite{MuV}, in
order to obtain all cyclic branched coverings of 2-bridge knots.

In Section~6 we will deal with the coloured graph representation
of manifolds introduced by M.~Pezzana \cite{Pe} and his school. A
wide class of coloured graph encoding 3-dimensional (possibly
singular) manifolds was defined by S.~Lins and A.~Mandel in
\cite{LM}, and intensively studied by many authors. It was shown
in \cite{Mu3} that all Lins-Mandel manifolds are cyclic branched
coverings of 2-bridge knots/links.

In Section~7 we will describe the presentation of the fundamental
group of these manifolds obtained in \cite{Mi} and \cite{Mu3}.
Moreover, in the case of coverings of 2-bridge knots, we also
describe the two different  cyclic presentations obtained in
\cite{Mi} and \cite{MuV}.

In Section~8 we give the first integer homology groups for a large
class of cases.

In Section~9 we will prove that each singly-cyclic branched
covering of a 2-bridge link is the composition of a
meridian-cyclic branched covering of a determined link and a
standard cyclic covering of $\S^3$ branched over a trivial knot.
The particular case of singly-cyclic coverings of the Whitehead
link was studied in \cite{CP}.

\section{Cyclic branched coverings}

Cyclic branched coverings of 2-bridge knots/links form a very
important class of 3-manifolds, that is a natural generalization
of the class of lens spaces. Indeed $M_{2,1}(\a/\b)$ is the lens
space $L(\a,\b)$ and, when $\gcd(n,k)=1$, $M_{n,k}(2/1)$ is the
lens space $L(n,k)$ \cite{Ro}. This class appears to be fairly
rich, since it contains several interesting 3-manifolds, such as
the Poincar\'e homology sphere, the Seifert--Weber hyperbolic
dodecahedron space, the Euclidean Hantzsche--Wendt manifold, the
hyperbolic Fomenko--Matveev--Weeks manifold and also an infinite
family of Brieskorn manifolds.

More precisely, the Poincar\'e homology sphere \cite{KS} is the
$5$-fold cyclic branched covering of the trefoil knot $\bt(3,1)$
and the $3$-fold cyclic branched covering of $\bt(5,1)$, i.e.,
$M_5(3/1)\cong M_3(5/1)$; the Seifert--Weber hyperbolic
dodecahedron space \cite{ST} is the $5$-fold singly-cyclic
branched covering of the Whitehead link $\bt(8,3)$ defined by
$k_1=1$ and $k_2=2$ , i.e., $M_{5,2}(8/3)$ (its generalizations
$M_{n,k}(8/3)$ were discussed in \cite{CP, HKM1, Zi2}); the
Euclidean Hantzsche--Wendt manifold \cite{Zi1} is the $3$-fold
cyclic branched covering of the figure-eight knot $\bt(5,3)$,
i.e., $M_3(5/3)$ (its generalizations $M_n(5/3)$, known as the
Fibonacci manifolds, were studied in \cite{CS, HKM2, HLM1, KV, MR,
MeV0, MeV1, RV}); the hyperbolic Fomenko--Matveev--Weeks manifold
\cite{FM, HoW}, which is the hyperbolic 3-manifold with the
smallest known volume, is the $3$-fold cyclic branched covering of
the knot $\bt(7,3)$, i.e., $M_3(7/3)$ \cite{MeV2} (its
generalizations $M_n(7/3)$ were studied in \cite{BKM, K2KV}); the
Brieskorn manifold $M(n,\a,2)$ \cite{M} is the $n$-fold
strictly-cyclic branched coverings of the torus knots or links
$\bt(\a,1)$, i.e., $M_n(\a/1)$ (see also \cite{CHK1, CHK2, Si}).

Note that, as a consequence of the positive solution of the Smith
conjecture \cite{MB}, this family of manifolds contains no sphere,
since $\bt(\a,\b)$ is never a trivial knot.

Now we list some sufficient homeomorphism conditions.

\begin{proposition}
Let $\bt(\a,\b)$ be a 2-bridge link and $k,k'\in\Z_n - \{ 0 \}$.
Then:
\begin{description}
\item[(i)] If $kk'=1$, then $M_{n,k}(\a/\b)$ is homeomorphic to
$M_{n,k'}(\a/\b)$;
\item[(ii)] $M_{n,k}(\a/\b)$ is homeomorphic to $M_{n,k}(-\a/\b)$;
\item[(iii)] $M_{n,k}(\a/\b)$ is homeomorphic to
$M_{n,-k}(\a/(\b-\a))$;
\item[(iv)] if the links $\bt(\a,\b)$ and $\bt(\a',\b')$ are
equivalent, then
$M_{n,k}(\a/\b)$ is homeomorphic to $M_{n,k}(\a'/\b')$.
\end{description}
\end{proposition}

\begin{proof}
Let $m_1$ and $m_2$ be meridians corresponding to the components
of $\bt (\a, \b)$. (i) Let $f : M_{n,k}(\a/\b) \to \S^3$ and $f' :
M_{n,k'}(\a/\b) \to \S^3$ the corresponding cyclic branched
coverings of $\bt(\a,\b)$. If $\ot$ and $\ot'$ are the associated
monodromy maps, we have $\ot[m_1]=1$, $\ot[m_2]=k$, $\ot'[m_1]=1$
and $\ot'[m_2]=k'$. Therefore, $\ot'\circ\r_\#=k'\ot$ and the two
coverings are equivalent. (ii) The link $\bt(\a,-\b)$ is
equivalent to the mirror image of $\bt(\a,\b)$. Let $\f$ be an
orientation reversing homeomorphism of $\S^3$ sending $\bt(\a,\b)$
to $\bt(\a,-\b)$. Since $\f$ preserve both the orientations of the
two components of $\bt(\a,\b)$, we have $\f_\#[m_1]=-[m'_1]$ and
$\f_\#[m_2]=-[m'_2]$, where $m'_1$ and $m'_2$ are the generator
meridians associated to $\bt(\a,-\b)$. Therefore, $\ot'\circ\f_\#
= - \ot$ and the two coverings are equivalent. (iii) The link
$\bt(\a,\b-\a)$ is equivalent to the link $\bt(\a,\b)$ with the
opposite orientation in the second component. Let $\f$ be the
identity map on $\S^3$, then $\f_\#[m_1]=[m'_1]$ and
$\f_\#[m_2]=-[m'_2]$. Therefore, $\ot'\circ\f_\#=\ot$ and the two
coverings are equivalent. (iv) Let $\f$ be the homeomorphism of
$\S^3$ realizing the equivalence. Without loss of generality we
can assume that $\f(K_1)=K'_1$ and $\f(K_2)=K'_2$. Then
$\f_\#[m_1]=[m'_1]$ and $\f_\#[m_2]=[m'_2]$. Therefore,
$\ot'\circ\f_\#=\ot$ and the two coverings are equivalent.
\end{proof}


\begin{corollary}
Let $\bt(\a,\b)$ be an arbitrary 2-bridge link and $k\in\Z_n -
\{ 0 \}$ such that $\gcd(n,k)=1$, then:
\begin{description}
\item[(i)] $M_{n,k}(\a/\b)$ is homeomorphic to $M_{n,k^{\pm
1}}(\a/\b)$;
\item[(ii)] if $\b^2=\a\pm 1$, then $M_{n,k}(\a/\b)$ is homeomorphic
to
$M_{n,\pm k^{\pm 1}}(\a/\b)$.
\end{description}
\end{corollary}

\begin{proof}
(i) See item (i) of the previous proposition.  (ii) Assume
$\b^2=\a+1$. Since $\b$ is odd, we have
$\b(\b-\a)=\b^2-\b\a=\b^2-\a=1$. Therefore, $\bt(\a,\b)$ is
equivalent to $\bt(\a,\b-\a)$ and from the previous proposition we
obtain $M_{n,k}(\a/\b)\cong M_{n,k}(\a/(\b-\a))\cong
M_{n,-k}(\a/\b)$. Now, let $\b^2=\a-1$. Since $\b$ is odd, we have
$\b(\a-\b)=\b\a-\b^2=\a-\b^2=1$. Therefore, $\bt(\a,\b)$ is
equivalent to $\bt(\a,\a-\b)$ and from the previous proposition we
obtain $M_{n,k}(\a/\b)\cong M_{n,k}(\a/(\a-\b))\cong
M_{n,k}(\a/(\b-\a))\cong M_{n,-k}(\a/\b)$.
\end{proof}

\medskip

When our manifolds have hyperbolic geometric structure, a partial
converse
to the previous results can be given (see Theorem \ref{Theorem
iff} below). Recall that a 2-bridge knot or link $\bt(\a,\b)$ is
hyperbolic (i.e., the complement $\S^3 - \bt(\a,\b)$ has
hyperbolic structure) if and only if $\bt(\a,\b)$ is non-toroidal,
that is $\b\not\equiv\pm 1$ mod $\a$. Thus, the geometric
structure of $M_{n,k}(\a/\b)$ can be obtained in these cases from
W.~Thurston \cite{Th} and W.~Dunbar \cite{Dn} results. Moreover,
when the branching set is toroidal, then $M_n(\a/\b)$ turns out to
be the Brieskorn manifold $M(n,\a,2)$. Thus, we have the following
result:

\begin{proposition} \label{Proposition geom}
\begin{description}
\item[(i)] \cite{HLM2} If $\gcd(n,k)=1$ and $\b \not \equiv \pm 1$
mod $\a$,
then $M_{n,k}(\a/\b)$ (or $M_{n}(\a/\b)$ in the case of a knot)
is hyperbolic for (i) $\a=5$, $n\ge 4$ and (ii)
$\a\ne 5$, $n\ge 3$. Moreover, $M_3(5/2)$ is Euclidean and $M_2
(\a/\b)$ is spherical for all $\a,\b$.
\item[(ii)] \cite{M} If $\b \equiv \pm 1$ mod $\a$, then
$M_n(\a/\b)$ is a
spherical manifold for $1/n + 1/\a > 1/2$, a Nil-manifold for $1/n
+ 1/\a = 1/2$, and a $\widetilde {SL}(2,\bf{R})$-manifold for $1/n
+ 1/\a < 1/2$.
\end{description}
\end{proposition}

From Theorem 1 of \cite{Zi2} (including the note (a) of page 293)
and Theorem 4.1 of \cite{Sa} (see tables of page 184), the next
result holds.

\begin{theorem}
\label{Theorem 2}
\label{Theorem iff}
Let $\bt(\a,\b)$ be a 2-bridge hyperbolic link. If $\gcd(n,k)=1$,
then $M_{n,k'}(\a/\b)$ and $M_{n,k}(\a/\b)$ are
homeomorphic if and only if
\begin{description}
\item[(i)] $k'=k^{\pm 1}$, when $\b^2\ne \a\pm 1$;
\item[(ii)] $k'=\pm k^{\pm 1}$, when $\b^2=\a\pm1$.
\end{description}
\end{theorem}


Volumes and Chern--Simons invariants of hyperbolic cyclic branched
coverings of 2-bridge knots were obtained in \cite{HLM3}, where the
table for small values of parameters is given.

Now we will show that every $n$-fold strictly-cyclic covering of a
$2$-bridge knot or link can be obtained as the 2-fold branched
covering of a certain $n$-periodic knot or link \cite{Do, SH}. If
$\rho : \S^3 \to \S^3$ is the involution previously described,
which leaves $\bt(\a,\b)$ invariant and exchanges its bridges, we
denote the quotient $\bt (\a,\b) / \langle \rho \rangle$ by $\bt
(1,\a,\b )$, that is the trivial knot pictured in Figure~3. So
$(\S^3, \bt (\a,\b))$ is a 2-fold covering of $(\S^3, \bt (1,
\a,\b))$ branched over a trivial knot $B$ corresponding to the
axis of the involution $\rho$. This covering gives us a natural
way to construct a periodic generalization of 2-bridge
knots/links. Denote by $\bt(n,\a, \b)$ the preimage of
$\bt(1,\a,\b)$ under the n-fold covering of $\S^3$ branched over
$B$. This knot or link admits a natural $n$-bridge presentation
(see an example in Figure~3) and will be called the {\it
$n$-cyclic extension} of $\bt (\a,\b)$. In particular,
$\bt(2,\a,\b) = \bt(\a,\b)$ and $\bt(n,\a, 1)$ is the torus
knot/link of type $(\a,n)$.

\begin{figure}
 \begin{center}
 \includegraphics*[totalheight=16cm]{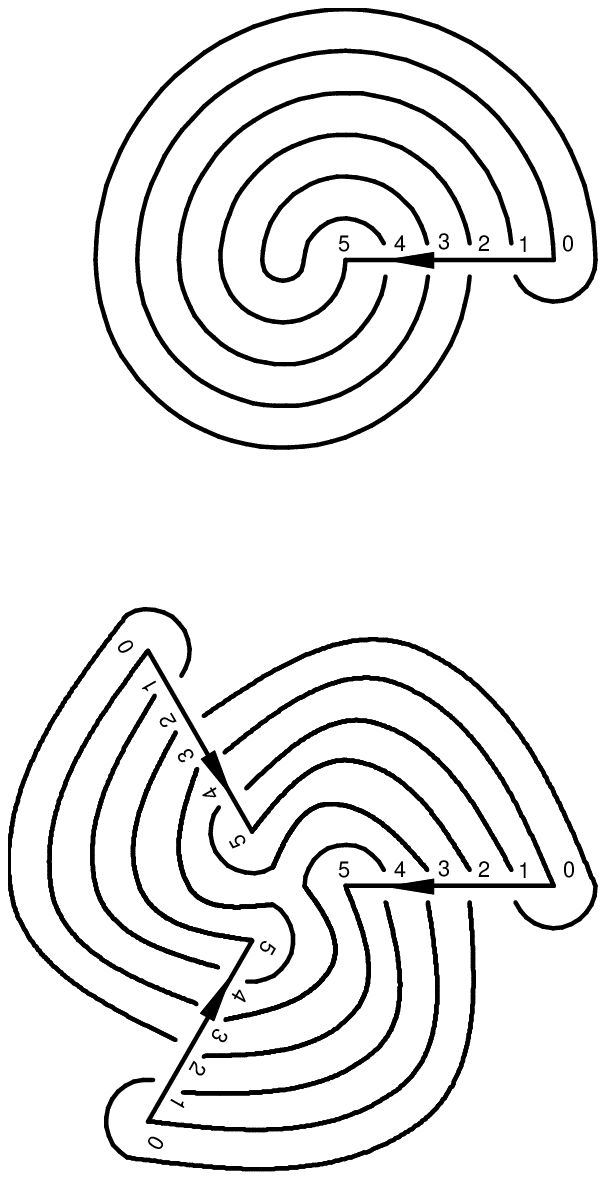}
 \end{center}
 \caption{Diagrams of $\bt(1, 5,3)$ and $\bt(3,5,3)$.}
 \label{Fig. 3}
\bigskip\bigskip
\end{figure}




\begin{proposition} \cite{Do, SH}
The manifold $M_n (\a / \b)$ is the $2$-fold branched co\-ve\-ring of
$\bt(n, \a, \b)$.
\end{proposition}

Thus, as an example, the Hantzsche--Wendt manifold is the $3$-fold
branched covering of $\bt (2,5,3) = \bt(5,3)$, pictured in
Figure~1, and the $2$-fold branched covering of
$\bt (3,5,3)$, pictured in Figure~3, that is a $3$-component
link known as the {\it Borromean rings} (see also \cite{MeV1}).

Independently, a $2$-fold covering description
of $M_n (\a, \b)$, when $\a$ is odd, arising from the generalized
Takahashi manifolds
(see Section~5), was obtained in \cite{MuV}.

\section{Polyhedral construction}

A classical method for constructing (orientable) closed $3$-manifolds
consists in the pairwise identification of (oppositely oriented)
boundary faces of a triangulated $3$-ball (see \cite{ST}).
The resulting quotient complex triangulates a closed pseudomanifold,
which is a  manifold if and only if its Euler characteristic
vanishes.

For strictly-cyclic branched coverings of 2-bridge knots/links
$\bt (p,q)$ this construction was realized by J. Minkus \cite{Mi}.
The idea is based on the method of realizing lens spaces. Here we
will briefly describe the construction.

Consider the boundary 2-sphere $\S^2$ of the 3-ball $B^3 = \{
(x,y,z) \in \R^3 \mid x^2 + y^2 +z^2 \le 1 \}$. Draw $n$ equally
spaced great semicircles joining the north pole $N=(0,0,1)$ to the
south pole $S=(0,0,-1)$. This decomposes $\S^2$ into $n$ congruent
lunes. Subdivide each semicircle into $p$ equal segments by
drawing $p -1$ equally spaced vertices on each semicircle. In
this way, each lune can then be viewed as a curvilinear $2p$
sided polygon on $\S^2$. Now bisect each lune by drawing a great
circle arc inside the lune, joining the vertex which is $q$
segments down from $N$ (the point $P_i$ of Figure~4) on each
semicircle with the vertex $q$ segments up from $S$ to the next
clockwise semicircle. Figure~4 shows this decomposition of
$\partial B^3=\S^2=\R^2\cup\{\infty\}$, where $S=\infty$.

\begin{figure}
 \begin{center}
 \includegraphics*[totalheight=10cm]{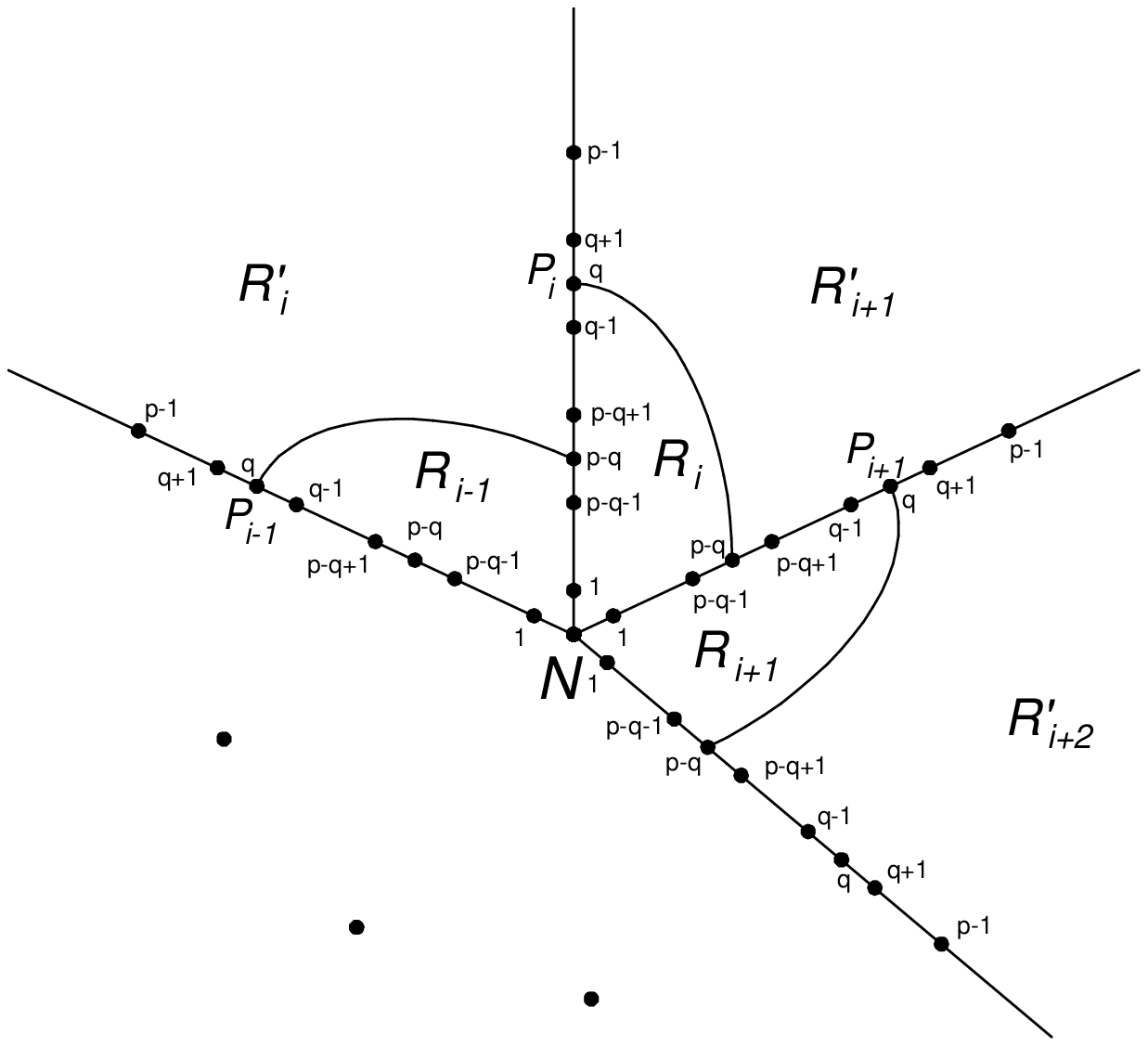}
 \end{center}
 \caption{Minkus polyhedral schemata for ${\widetilde M}_n(p,q)$.}
 \label{Fig. 4}
\bigskip\bigskip
\end{figure}




The result is the decomposition of $\partial B^3$ into $2n$
regions $R_i$, $R'_i$, $i=1,\ldots, n$. The regions $R_i$ are
around $N$ and the regions $R'_i$ are around $S$, and each $R'_i$
can be reached from the corresponding $R_i$ by moving $R_i$
counterclockwise to the adjacent lune and then shifting from the
northern to the southern hemisphere of $\partial B^3$. The cell
3-complex ${\widetilde M}_n(p,q)$ is obtained from $B^3$ by
identifying $R_i$ with $R'_i$ on $\partial B^3$ for each $i=1,
\ldots , n$  by an orientation reversing homeomorphism which
matches the vertex $P_i$ of $R_i$ with the vertex $P_{i-1}$ of
$R'_i$.

\begin{theorem} \label{Theorem Minkus} \cite[Theorem~7]{Mi}
The manifold ${\widetilde M}_n (p,q)$ is the $n$-fold
strictly-cyclic branched covering of the 2-bridge knot or link
$\bt(p,q)$.
\end{theorem}

Therefore the manifolds ${\widetilde M}_n (p,q)$ are homeomorphic
to the above defined manifolds $M_n(p/q)$ and provide a polyhedral
construction for them.
As an example, Figure~5 gives the polyhedral construction of
the Hantzsche--Wendt manifold, that is $M_3(5/3)$ in our notations.

\begin{figure}
 \begin{center}
 \includegraphics*[totalheight=6cm]{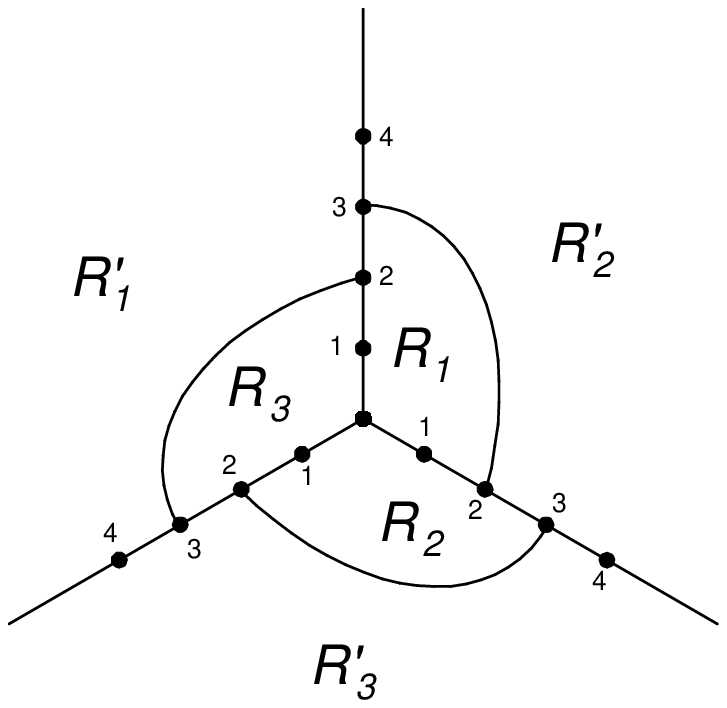}
 \end{center}
 \caption{Minkus polyhedral schemata for the Hantzsche--Wendt
manifold.}
 \label{Fig. 5}
\bigskip\bigskip
\end{figure}




The Minkus polyhedral construction only produces strictly-cyclic
branched coverings of 2-bridge knots/links. The generalization to
singly-cyclic coverings $M_{n,k}(p/q)$ is straightforward and is
illustrated in Figure~6, which is the same as Figure 5 of \cite{Mu3}.

In this case, the cell 3-complex $M_{n,k}(p/q)$ is obtained
from $B^3$ by identifying $R_i$ with $R'_i$ on $\partial B^3$
for each $i=1, \ldots , n$  by an orientation reversing
homeomorphism which matches the vertex $P_i$ of $R_i$ with
the vertex $P_{i-k}$ of $R'_i$.

\begin{figure}
 \begin{center}
 \includegraphics*[totalheight=10cm]{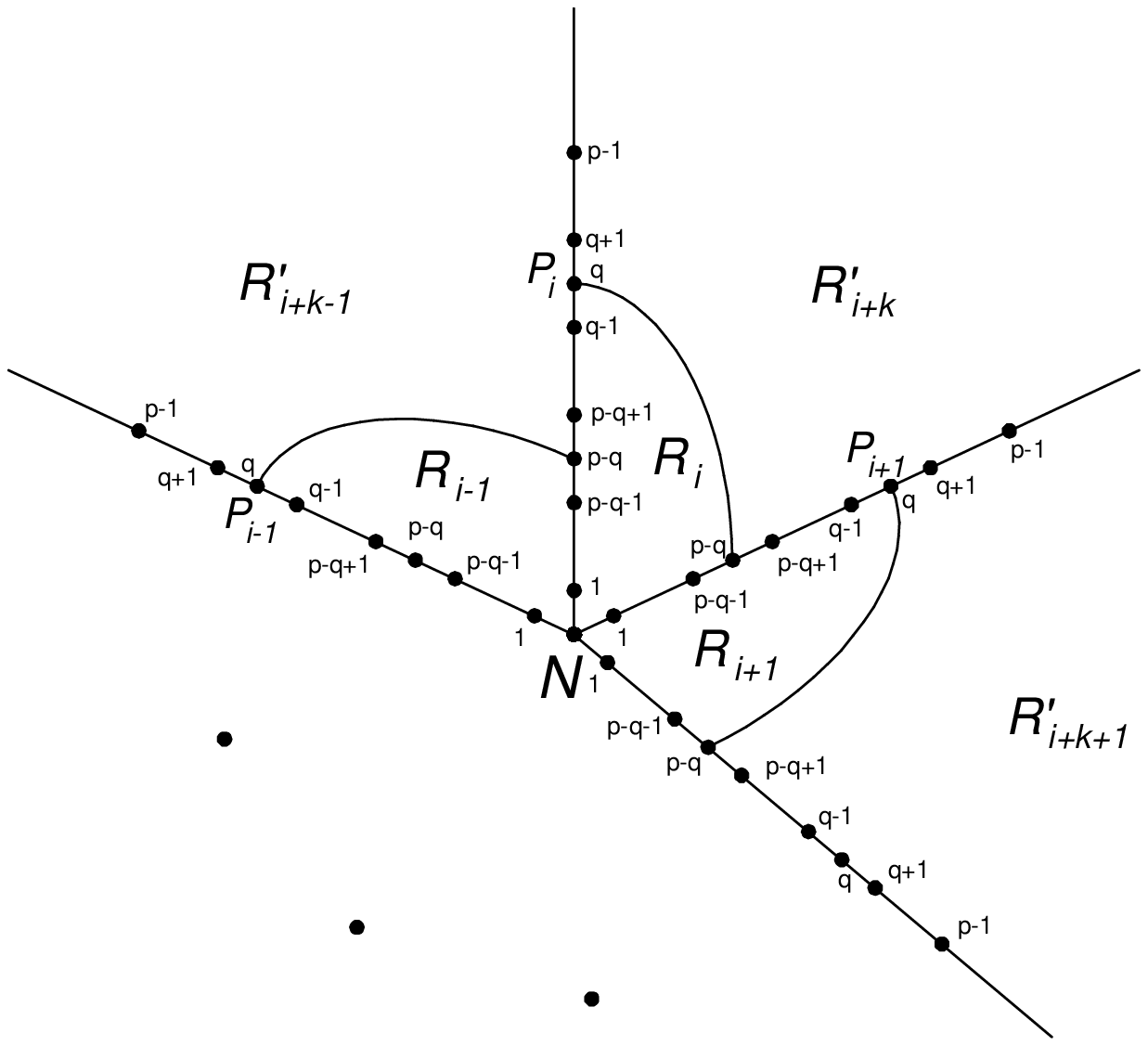}
 \end{center}
 \caption{Polyhedral schemata for $M_{n,k}(p/q)$.}
 \label{Fig. 6}
\bigskip\bigskip
\end{figure}




We remark that another polyhedral construction for
$M_n((2\ell-1)/\ell)$ is presented in \cite{KiY}.

\section{Heegaard diagram construction and genus}

In this section we will discuss Heegaard diagrams
as well as cyclically symmetric Heegaard diagrams for
our manifolds.

Heegaard diagrams for meridian-cyclic branched coverings of
2-bridge knots/links can be directly obtained from the coloured
graph construction that will be discussed in Section~6 (see
\cite{Mu3}). More precisely, a Heegaard diagram of genus $n-1$ of
$M_{n,k}(\a / \b)$, with $\gcd(n,k)=1$, can be obtained from
Figure~7 by identifying the disk $C_i$ with the disk $C'_i$, for
$i=1, \ldots, n-1$, according to the numeration of the vertices,
and removing one of the $n$ closed curves arising from these
identifications, illustrated by the dashed lines in Figure~7.

\begin{figure}
 \begin{center}
 \includegraphics*[totalheight=18cm]{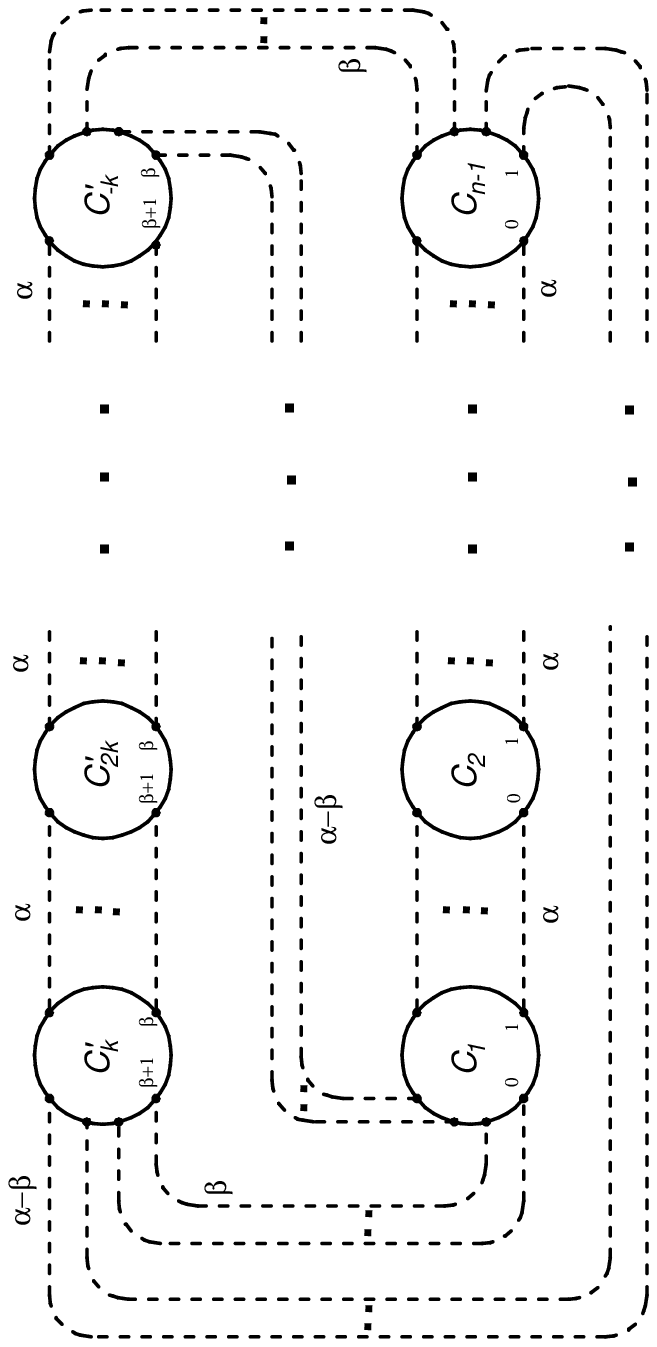}
 \end{center}
 \caption{}
 \label{Fig. 7}
\bigskip\bigskip
\end{figure}




Figure~8 presents a Heegaard diagram of $M_3(5/3)$
obtained in this way.

\begin{figure}
 \begin{center}
 \includegraphics*[totalheight=8cm]{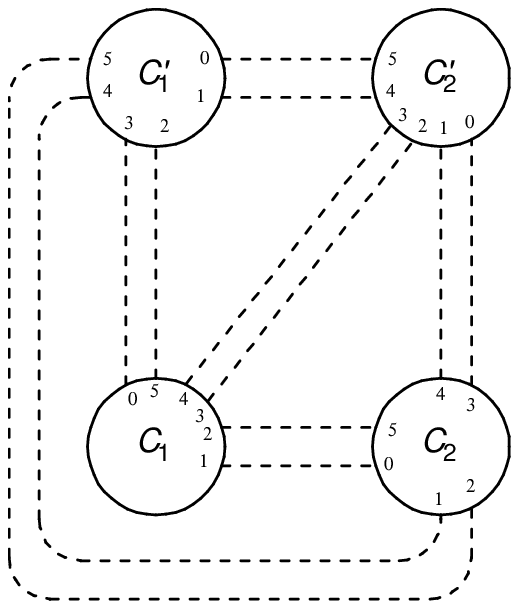}
 \end{center}
 \caption{An Heegaard diagram of genus 2 for the Hantzsche--Wendt
manifold.}
 \label{Fig. 8}
\bigskip\bigskip
\end{figure}




At the same time, it is natural to ask about Heegaard diagrams for
our manifolds corresponding to cyclic coverings. The connections
between cyclic branched coverings of knots and cyclic
presentations of groups induced by suitable Heegaard diagrams have
recently been discussed in several papers
\cite{BKM,CHK1,CHK2,CHR,Du,HKM1,HKM2,Ki,KV,K2KV,MR}. In order to
investigate these relations, M.J. Dunwoody introduced in \cite{Du}
a class of Heegaard diagrams, depending on six integers, having a
cyclic symmetry and encoding a cyclic presentation for the
fundamental group of the represented manifold $D(a,b,c,n,r,s)$. In
\cite{GM} it has been shown that the 3-manifolds represented by
these diagrams (called {\it Dunwoody manifolds}) are cyclic
coverings of lens spaces branched over genus one 1-knots (also
called $(1,1)$-knots). As a corollary, it has been demonstrated
that, for particular values of the parameters, the Dunwoody
manifolds turn out to be cyclic coverings of ${\bf S}^3$ branched
over some knots. This gives a positive answer to a conjecture made
by Dunwoody, which has also been independently proved in
\cite{SK}.

It is interesting to note that the class of Dunwoody manifolds
properly contains the class of cyclic branched coverings of
2-bridge knots, as stated in the following theorem. Note that in
the statement $\bar s$ is an integer only depending on $a$ and $r$
(see details in Section 3 of \cite{GM}).

\begin{theorem} \label{Theorem Dunwoody} \cite{GM}
For all $a,r>0$ and $n>1$, the $n$-fold cyclic covering of $\S^3$
branched over $\bt (2a+1,2r)$ is the Dunwoody manifold
$D(a,0,1,n,r,\bar s)$. Thus, all branched cyclic coverings of
2-bridge knots are Dunwoody manifolds.
\end{theorem}

\begin{figure}
 \begin{center}
 \includegraphics*[totalheight=8cm]{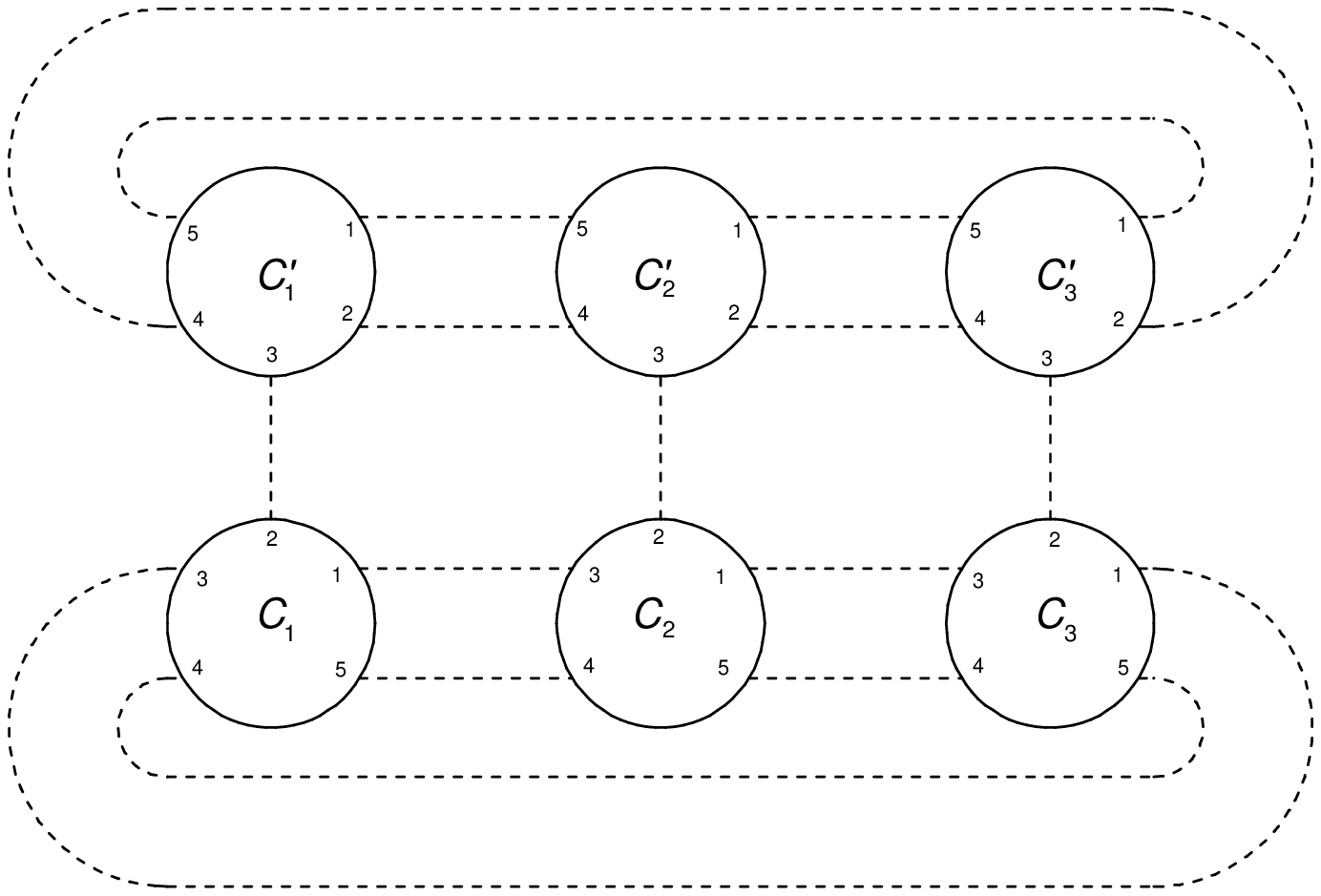}
 \end{center}
 \caption{A 3-symmetric Heegaard diagram of genus 3 for $D(2,0,1,3,1,0)$, the Hantzsche--Wendt
manifold.}
 \label{Fig. 9}
\bigskip\bigskip
\end{figure}




As a consequence of the previous theorem, each $n$-fold cyclic
branched covering of a 2-bridge knot admits a Heegaard
diagram of genus $n$ with a cyclic symmetry of order $n$.
These diagrams induce ``geometric'' cyclic presentations
for the corresponding fundamental groups (see Section~7).

Now we discuss some results on genus (classical, $p$-symmetric
and equivariant) of our manifolds.

{}From the standard representation of a 2-bridge knot/link as a
4-plat $(\S^3,\bt(\a,\b))=(B',A')\cup_f(B'',A'')$, where $B'$,
$B''$ are 3-balls and $A'\subset B'$, $A''\subset B''$ are pairs
of trivially properly embedded arcs, we get a Heegaard splitting
of $M_{n,k',k''}(\a/\b)$, where the splitting surface is an
$n$-fold cyclic covering of $\S^2$ branched over 4 points. From
the Riemann-Hurwitz formula, the surface has genus
$g=n+1-\gcd(n,k') - \gcd(n,k'')$ (see details in \cite{BZ}).
Therefore we have the following:

\begin{proposition} \label{Proposition genus}
Let $g(M)$ be the genus of a 3-manifold $M$. Then
$$
g(M_{n,k',k''}(\a/\b))\le n+1- \gcd(n,k') - \gcd(n,k'') .
$$
In particular, for singly-cyclic coverings
$ g(M_{n,k}(\a/\b))\le n- \gcd(n,k) $ and for meridian-cyclic
coverings $g(M_{n,k}(\a/\b))\le n-1$. Thus, for strictly-cyclic
coverings $g(M_n(\a/\b))\le n-1$.
\end{proposition}

\medskip

{}From Theorem 4 of \cite{Mu4}, for the case of strictly-cyclic
branched coverings the same estimation also holds for the
$n$-symmetric genus introduced by J.~Birman and H.~Hilden in
\cite{BH}.

\begin{proposition} \label{Proposition symm-genus}
Let $g_n(M)$ be the $n$-symmetric genus of a 3-manifold $M$. Then
$$g_n(M_{n}(\a/\b))\le n-1.$$
\end{proposition}

In many interesting cases the result of Proposition
\ref{Proposition genus} can be improved.

\begin{proposition} \label{Proposition braid}
\begin{description}
\item[(i)] For all $n,\a>1$ we have
$$g(M_n(\a/1))\le\min\{\a-1,n-1\}.$$
\item[(ii)] For all $n,c>1$ we have
$$g(M_n((3c-1)/3))\le \min\{c,n-1\}.$$
\end{description}
\end{proposition}

\begin{proof}
(i) The $n$-fold strictly-cyclic co\-vering of $\bt(\a,1)$ is the
2-fold covering of the torus knot/link of type $(\a, n)$ \cite{M},
which can be presented as the closure of a $\a$-string braid. So
it has bridge number $b \le \a$ and then $g \le \a -
1$ \cite{BH}. (ii) The $n$-fold strictly-cyclic covering of
$\bt(3c-1,3)$ is the 2-fold covering of $\bt(n,3c-1,3)$, which can
be presented as the closure of a $(c+1)$-string braid. So
it has bridge number $b\le c+1$ and then $g \le c$
\cite{BH}.
\end{proof}

\medskip

In particular, from (ii) we see that the $n$-fold cyclic branched
covering of the figure-eight knot ($c=2$) has genus $g=2$ if $n>2$
and $g=1$ if $n=2$ (see also \cite{MeV3}).
Similarly, the $n$-fold strictly-cyclic branched
coverings of the Whitehead link ($c=3$) have genus $g\le 3$.

\smallskip

Recall that the maximally possible order of a finite group
of orientation-preserving homeomorphisms of a handlebody
$V_g$ of genus $g \ge 2$ equals $12(g-1)$ (see \cite{Zi5}),
analogously to the classical Hurwitz $84(g-1)$-bound for
a closed Riemann surfaces of genus $g\ge 2$.

According to \cite{Zi3}, a closed 3-manifold $M$ is said to be a
{\it $G$-manifold of genus $g$} if it admits an action of the
finite group $G$ and $g$ is the minimal genus of a Heegaard
splitting of $M$ for which both handlebodies are invariant under
the $G$-action. A $G$-manifold of genus $g>1$ is called {\it
minimal} if the induced $G$-action on each of the two handlebodies
of an invariant Heegaard splitting of genus $g$ is a strong genus
action (i.e., there is no action of $G$ on a handlebody of genus
$\bar g$, with $1<\bar g<g$). Moreover, if $G$ has maximal
positive order $12(g-1)$ then the $G$-manifold $M$ and the
$G$-action are called {\it maximally symmetric}. For $n>2$, every
minimal $\Z_n$-manifold is a minimal $\Di_n$-manifold, where
$\Di_n$ is the dihedral group of order $2n$. If $n>2$ is prime,
every minimal $\Z_n$-manifold is a minimal
$(\Di_n\times\Z_2)$-manifold \cite{Zi3}. Each manifold of genus 2
is a $G$-manifold of genus 2, where $G$ is one of the four groups
$\Z_2$, $\Di_2$, $\Di_4$ or $\Di_6$, and the $\Di_6$-manifolds are
maximally symmetric \cite{Zi3}.

The cyclic branched coverings of 2-bridge knots/links play
an important role in this theory.

\begin{theorem} \label{Theorem minimal} \cite{Zi3}
\begin{description}
\item[(i)] The minimal $(\Di_n\times\Z_2)$-manifolds, of genus
$g=n-1$,
are exactly the $n$-fold strictly-cyclic branched coverings of
2-bridge knots/links. The minimal $\Z_n$ or $\Di_n$-manifolds that
are not minimal $(\Di_n\times\Z_2)$-manifolds are cyclic
branched coverings of 2-bridge links with two components of
different branching index.
\item[(ii)] The maximally symmetric $\Di_6$-manifolds are
exactly the
3-fold cyclic branched coverings of 2-bridge knots/links.
\end{description}
\end{theorem}

\section{Surgery construction}

In this section we describe a symmetric surgery presentation of
cyclic branched coverings of 2-bridge knots. This presentation
comes from a more general construction of {\it generalized
Takahashi manifolds} introduced in \cite{MuV}. For any pair of
positive integers $m$ and $n$, we consider the link ${\cal
L}_{n,m} \subset \S^3$ with $2mn$ components drawn in Figure~10.
All its components $c_{i,j}$, $1\le i \le 2n$, $1 \le j \le m$,
are unknotted circles and they form $2n$ subfamilies of $m$
unlinked circles $c_{i,j}$, $1 \le j \le m$, with a common center.
The link ${\cal L}_{n,m}$ has a cyclic symmetry of order $n$ which
permutes these $2n$ subfamilies of circles.

Consider the manifold obtained by Dehn surgery on $\S^3$, along
the link ${\cal L}_{n,m}$, such that the surgery coefficients
$p_{k,j} / q_{k,j}$ correspond to the components $c_{2k-1,j}$, and
$r_{k,j} / s_{k,j}$ correspond to the components $c_{2k,j}$, where
$1 \le k \le n$ and $1 \le j \le m$.
Without loss of generality, we can always suppose that
$\gcd(p_{k,j},q_{k,j})=1$, $\gcd(r_{k,j},s_{k,j})=1$ and
$p_{k,j},r_{k,j}\ge 0$.

We will denote the resulting 3-manifold by $T_{n,m} (p_{k,j} /
q_{k,j} ; r_{k,j} / s_{k,j} )$. This manifold will be called
a {\it generalized Takahashi manifold} since for $m=1$ we get the
manifolds introduced by M.~Takahashi in \cite{Ta}.

When the surgery coefficients are $n$-periodic, i.e.
$p_{k,j}=p_j$, $q_{k,j}=q_j$, $r_{k,j}=r_j$, and $s_{k,j}=s_j$,
the resulting manifold $T_{n,m} (p_1/q_1,\ldots,p_m/q_m;
r_1/s_1,\ldots,r_m/s_m)$ is said
to be a {\it generalized periodic ($n$-periodic) Takahashi
manifold}.

The following theorem generalizes the result obtained in \cite{Ta}
for periodic Takahashi manifolds, and gives the relation between
the generalized periodic Takahashi manifolds and the cyclic
branched coverings of 2-bridge knots.

\begin{theorem} \label{twobridge}  \cite{MuV}
The manifold $T_{n,m}(1/q_1,\ldots,1/q_m; 1/s_1,\ldots,1/s_m)$ is
homeomorphic to $M_n (\a/\b)$, where $\a/\b$ is the rational
number defined by the continued fraction
$[-2q_1,2s_1,\ldots,-2q_m,2s_m]$.
\end{theorem}

Because every 2-bridge knot admits a Conway representation with an
even number of even parameters (see Exercise 2.1.14 of
\cite{Ka}), we have the following property:

\begin{corollary} \label{all twobridge}  \cite{MuV}
The family of generalized periodic Takahashi manifolds contains
all the cyclic branched coverings of 2-bridge knots.
\end{corollary}

As a result, a surgery presentation with cyclic symmetry is
obtained for all cyclic branched coverings of 2-bridge knots (see
Figure~10).

\begin{figure}
 \begin{center}
 \includegraphics*[totalheight=7cm]{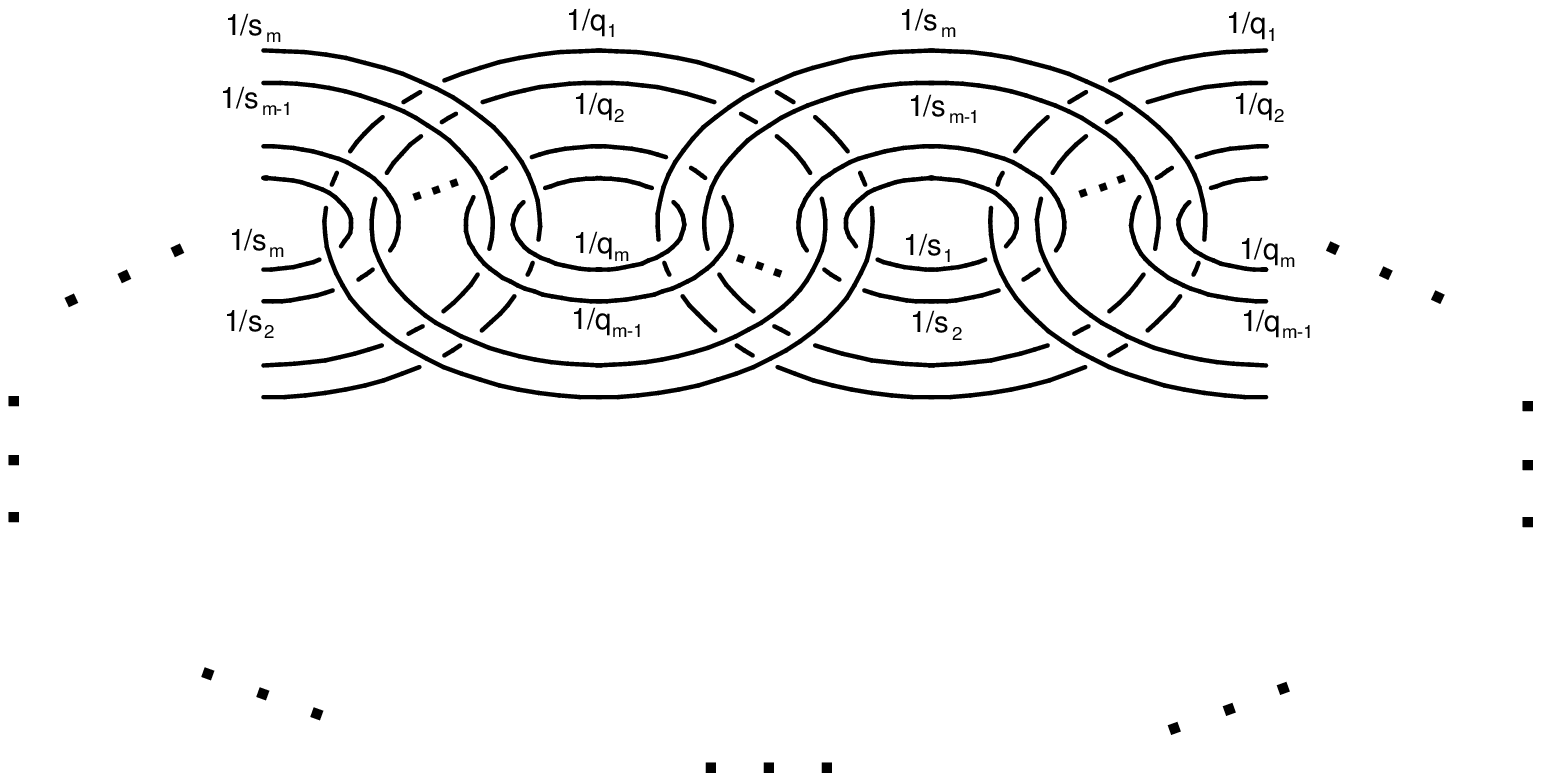}
 \end{center}
 \caption{Surgery presentation of $M_n(\a/\b)$, where
$\a/\b=[-2q_1,2s_1,\ldots,-2q_m,2s_m])$.}
 \label{Fig. 10}
\bigskip\bigskip
\end{figure}




As an example, the surgery description of the Hantzsche--Wendt
manifold (that is $M_3(5/2)$ in the above notation) is illustrated
in Figure~11.

\bigskip\bigskip

\begin{figure}
 \begin{center}
 \includegraphics*[totalheight=5cm]{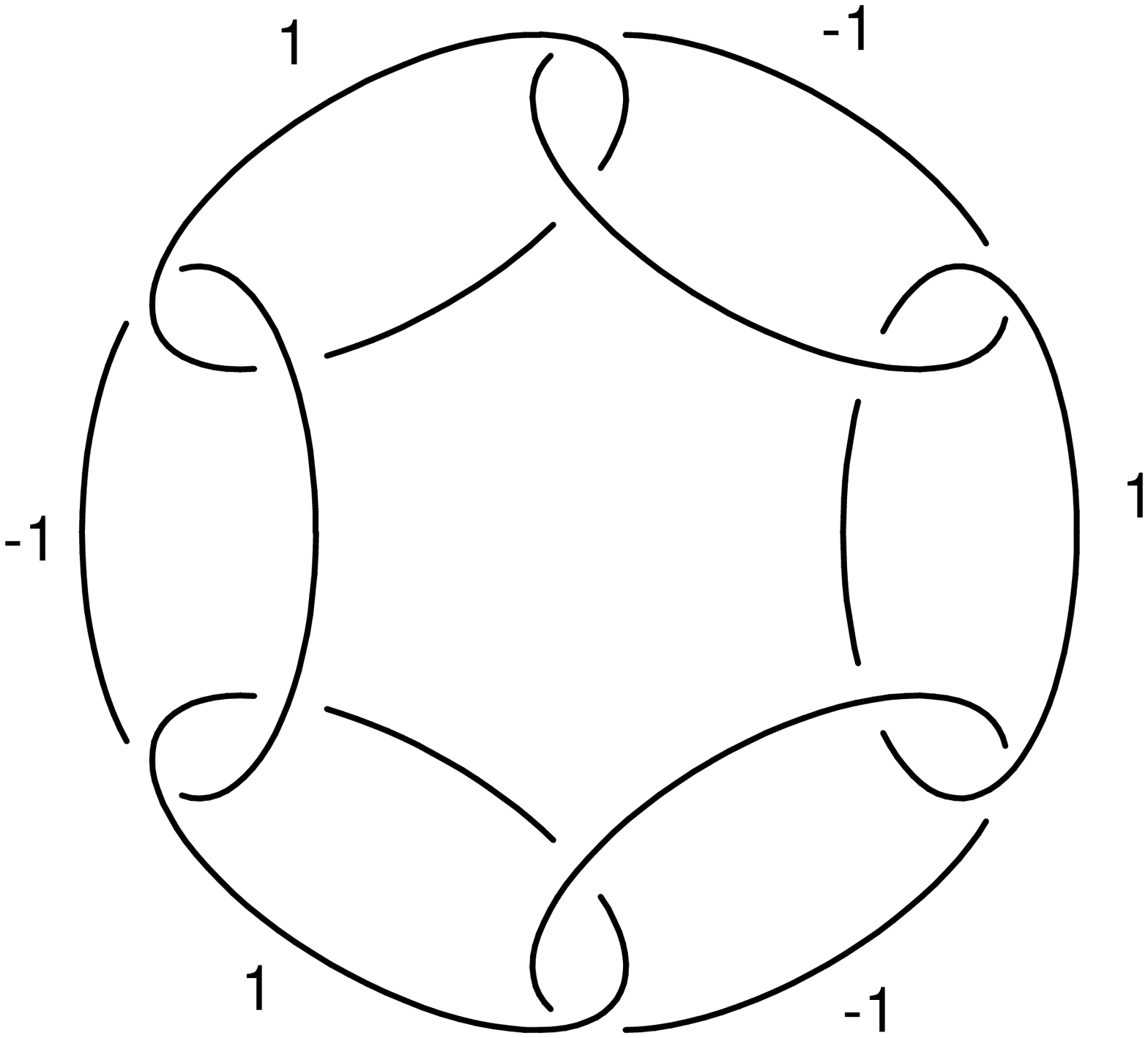}
 \end{center}
 \caption{Surgery description of the Hantzsche--Wendt manifold.}
 \label{Fig. 11}
\bigskip\bigskip
\end{figure}




\section{Coloured graph construction}

We recall some facts about the coloured-graph presentation of
3-manifolds.

A $4$-{\it coloured graph} is a pair $(\Gamma ,\gamma )$, where
$\Gamma$ is a finite connected $4$-regular graph without loops,
and $\gamma :E(\Gamma )\to\Delta_3=\{0,1,2,3\}$ is a proper
edge-coloration (i.e., adjacent edges have different colours).
Every $4$-coloured graph represents a pseudosimplicial complex
\cite{HW} $K(\Gamma )$ defined in the following way: (i) take a
$3$-simplex $\sigma (x)$ for each vertex $x\in V(\Gamma )$ and
label its vertices by the elements of $\Delta_3$; (ii) identify,
for every pair $x,y\in V(\Gamma )$ of $c$-adjacent vertices, the
$2$-faces of $\sigma (x)$ and $\sigma (y)$ opposite the
vertices labelled by $c$. The underlying space $\vert K(\Gamma
)\vert$ is a connected $3$-dimensional pseudomanifold, which is
orientable if and only if $\Gamma$ is bipartite \cite{FGG}. Note
that only isolated singular points may appear, and such spaces are
said to be {\it singular manifolds} in \cite{Mo2}. A $3$-{\it
gem\/} is a $4$-coloured graph representing a $3$-manifold; every
manifold $M$ is representable by gems \cite{FGG, Li, Pe}. A gem is
called a {\it crystallization} if, for each colour $c \in
\Delta_3$, the subgraph of $\Gamma$ obtained by removing all
$c$-coloured edges is connected.

There is a strict connection between crystallizations and Heegaard
diagrams of 3-manifolds \cite{Pe}. Let $c'$  and $c''$ be two
colours from $\Delta_3$, and $A=\{ c', c'' \}$, $B = \Delta_3 -
A$. A Heegaard diagram of the represented manifold can be obtained
just by removing one $A$-cycle (i.e., formed by edges coloured by
the colours of $A$) and one $B$-cycle. The remaining $A$-cycles
and $B$-cycles give, respectively, the two systems of curves of
the Heegaard diagram. In other words, the $A$-cycles and the
$B$-cycles of a crystallization represent an extended Heegaard
diagram of the manifold (see \cite{NO}).

Also manifolds of arbitrary dimension can be represented by
edge-coloured graphs, which give a combinatorial way of
representing manifolds; for general references see \cite{FGG} and
\cite{Li}.

\smallskip

The family of Lins-Mandel 4-coloured graphs $G(n,p,q,c)$, with
$n,p>0$, $q\in \Z_{ 2p}$, $c\in \Z_n$ and $\gcd(p,q)=1$, has been
defined in \cite{LM}. The set of vertices of $G(n,p,q,c)$ is
$V=\Z_n \times \Z_{2p}$ and the coloured edges are obtained by the
following four fixed-point-free involutions on $V$: $$
\begin{array}{ll}
\iota_0(i,j)=(i+c\eta(j-q),1-j+2q), & \qquad \iota_1(i,j) =(i+\eta
(j),1-j),\\ \iota_2(i,j)=(i,j+(-1)^j), & \qquad
\iota_3(i,j)=(i,j-(-1)^j);
\end{array}
$$ where $\eta : \Z_{2p} \to \{-1,+1\}$ is the map defined by $$
\eta (j)= \left\{
\begin{array}{ll}
+1  & \mbox{if }1 \le j \le p  \\
-1 & \mbox{otherwise}
\end{array} \right. \, .
$$ For each $k\in\{0,1,2,3\}$, we join the vertices $v,w\in V$ by
a $k$-coloured edge if and only if $\iota_k(v)=w.$

Roughly speaking, the graph $G(n,p,q,c)$ is obtained by taking $n$
copies $C_i$, $i\in\Z_n$, of the $\{2,3\}$-cycle of length $2p$
(so that $V(C_i)=\{(i,j)\mid j\in \Z_{2p}\}$) joined with
$C_{i-1}$ and $C_{i+1}$ by $p$ edges of colour 1, and with
$C_{i-c}$ and $C_{i+c}$ by $p$ edges of colour 0.

Each graph $G(n,p,q,c)$ is connected and bipartite; hence, it
represents a connected, orientable 3-dimensional pseudomanifold
$S(n,p,q,c)$. This class of graphs and spaces have been
intensively studied \cite{CG1, CG2, Ca1, Ca2, Ca3, Ca4, Do, Gr,
JT, Li, LM, LMu, Mu2, Mu3, Mu1}. Remark that Lins-Mandel spaces
have been introduced as a combinatorial generalization of the lens
spaces, since $G(2,p,q,1)$ is the standard graph representing the
lens space $L(p,q)$ (see \cite{DG}).

\begin{figure}
 \begin{center}
 \includegraphics*[totalheight=15cm]{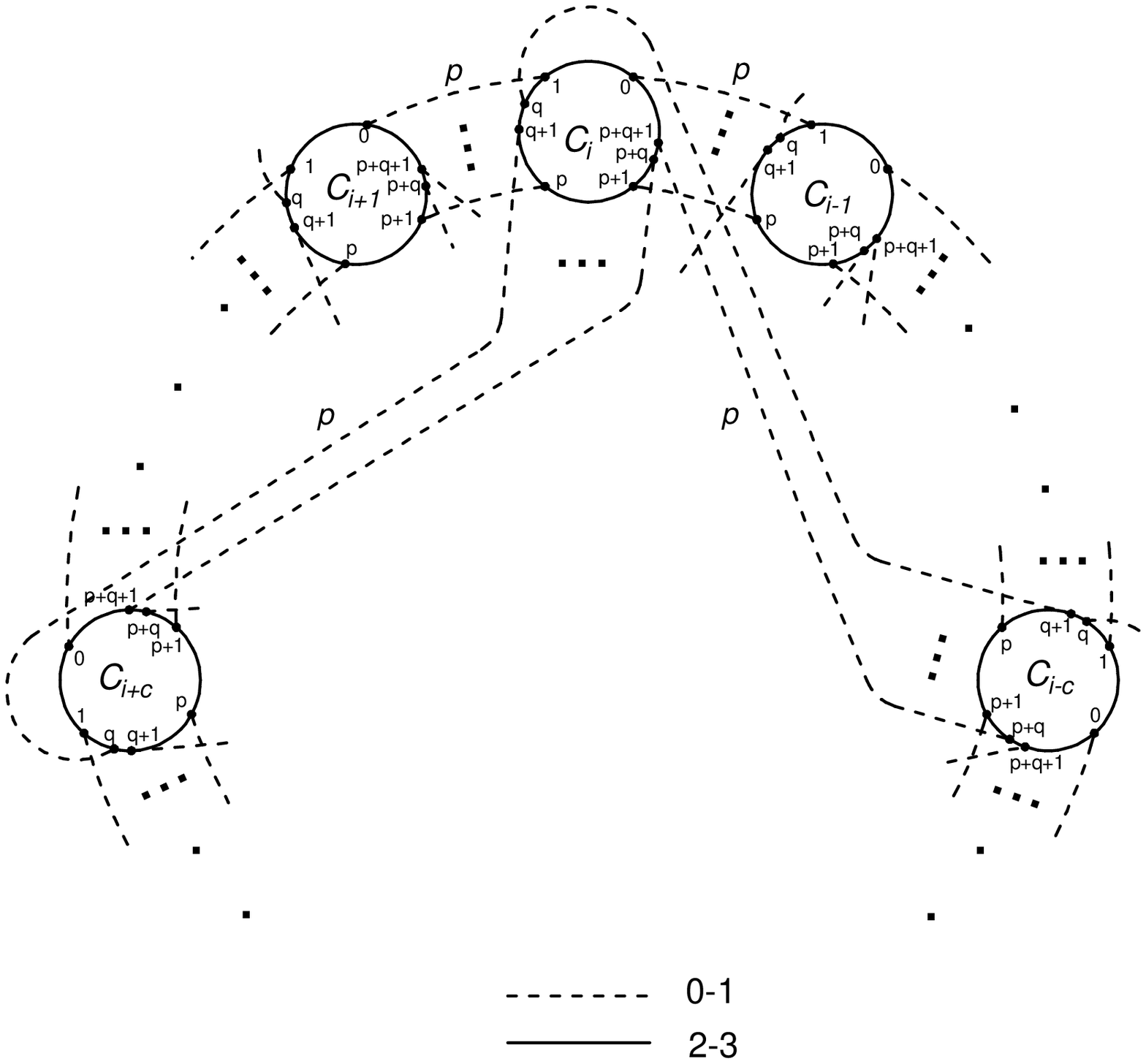}
 \end{center}
 \caption{The Lins--Mandel graph $G(n,p,q,c)$.}
 \label{Fig. 12}
\bigskip\bigskip
\end{figure}

The following characterization of the Lins-Mandel gems exists:

\begin{proposition} \cite{Mu2}
A Lins-Mandel graph $G(n,p,q,c)$ represents a 3-manifold if and
only if either $p$ is even or $c=0,(-1)^q$.
\end{proposition}

Moreover, a Lins-Mandel gem $G(n,p,q,c)$ is a crystallization
if and only if $\gcd(n,c) = 1$ \cite{CG1}.

\begin{theorem} \label{Lins-Mandel as cyclic} \cite{Mu3, Mu1}
Lins-Mandel spaces have the following
topological structure:
\begin{description}
\item[(i)] If $S(n,p,q,c)$ is a manifold with $c\ne 0$ and $p>1$,
then
it is homeomorphic to $M_{n,-c}(p/q)$, otherwise
it is homeomorphic to $\S^3$.
\item[(ii)] If $S(n,p,q,c)$ is not a manifold, then it is an $n$-fold
cyclic covering of $\S^3$ branched over a $\theta$-graph, which is
embedded as the  2-bridge knot $\bt(p,q)$ with an unknotting tunnel.
\end{description}
\end{theorem}

Thus, the manifold $S(n,p,q,c)$ is a singly-cyclic branched
covering of $\bt (p,q)$. In particular the covering is
strictly-cyclic when $c=(-1)^q$, is almost-strictly-cyclic when
$c=\pm 1$, and is meridian-cyclic when $\gcd(n,c)=1$.

Remark that by removing the cycle $C_n$ from the graph
$G(n,p,q,c)$, with $\gcd(n,c)=1$, we will get exactly the graph
pictured in Figure~7, once the relative identifications of $C_i$
with $C'_i$, for $i=1, \ldots, n-1$, have been performed,
according to the numeration of the vertices.

Necessary and sufficient conditions for the isomorphism between
Lins-Mandel gems are obtained in \cite{LMu}, when $n,p>2$. As a
consequence, sufficient conditions for the homeomorphism between
Lins-Mandel manifolds not homeomorphic to the sphere or the lens
spaces directly follow.

\begin{theorem} \cite{LMu} \label{Theorem 1}
For $n,p > 2 $ the following isomorphisms of graphs hold:
\begin{description}
\item[(i')] If $p$ is even and $\gcd(n,c)\ne 1$,
then $G(n',p',q',c')\cong G(n,p,q,c)$ if and only if
$$\mbox{ (1) } n'=n,
\mbox{ (2) } p'=p, \mbox { (3)  either } \left\{\begin{array}{lc}
q'=\pm q^{\pm 1}\\c'=c \\
\end{array}\right.\mbox  { or }
\left\{\begin{array}{lc} q'=\pm q^{\pm 1}+p\\c'=-c\\
\end{array}\right..$$
\item[(i'')] If $p$ is even and $\gcd(n,c)=1$, then
$G(n',p',q',c')\cong G(n,p,q,c)$ if and only if
$$\mbox{ (1) } n'=n,
\mbox{ (2) } p'=p, \mbox { (3)  either }\left\{\begin{array}{lc}
q'=\pm
q^{\pm 1}\\c'=c^{\pm 1}\\
\end{array}\right. \mbox  { or }
\left\{\begin{array}{lc} q'=\pm q^{\pm 1}+p\\c'=-c^{\pm 1}\\
\end{array}\right..$$
\item[(ii)] If $p$ is odd, then
$G(n',p',q',(-1)^{q'})\cong G(n,p,q,(-1)^q)$ if and only if
$$\mbox{ (1) } n'=n, \mbox{ (2) } p'=p, \mbox{ (3) }
q'\equiv\pm q^{\pm 1} \mbox{ mod } p.$$
\end{description}
Hence, if one of the above conditions holds, then the manifolds
$S(n,p,q,c)$ and $S(n',p',q',c')$ are homeomorphic.
\end{theorem}

Observe that the isomorphism conditions of part (ii) of the previous
theorem
are the same as the homeomorphism conditions for
lens spaces. This is not true for part (i), since, in this case, the
situation is complicated by the presence of the additional parameter
$c$.

Cases where $p$ is even are particularly interesting because the
graph
always represents a manifold without any restriction on $c$. From
Theorem \ref{Theorem 1} we get:

\begin{corollary} \cite{LMu} \label{Corollary 1}
Let $n,p,q$ be fixed,  with $n,p\ge 3$ and $p$ even.
Then $G(n,p,q,c')\cong G(n,p,q,c)$ if and only if
\begin{description}
\item[(i)] $c'=c$, $\,$ when $\gcd(n,c)\ne 1$ and $q^2\ne p\pm 1$;
\item[(ii)] $c'=\pm c$, $\,$ when $\gcd(n,c)\ne 1$ and $q^2=p\pm 1$;
\item[(iii)] $c'=c^{\pm 1}$, $\,$ when $\gcd(n,c)=1$ and $q^2\ne
p\pm 1$;
\item[(iv)] $c'=\pm c^{\pm 1}$, $\,$ when $\gcd(n,c)=1$ and
$q^2=p\pm 1$.
\end{description}
\end{corollary}

Due to Theorem \ref{Theorem 2} and Corollary \ref{Corollary 1}, in
many cases graphs distinguish manifolds:

\begin{corollary} \label{Theorem 3}
Assume $n,p > 2$ and $q\ne\pm 1,p\pm 1$.
For $\gcd(n,c)=1$, the manifolds $S(n,p,q,c')$ and $S(n,p,q,c)$
are homeomorphic if and only if the graphs $G(n,p,q,c')$ and
$G(n,p,q,c)$ are isomorphic.
\end{corollary}

The Lins-Mandel family contains only singly-cyclic coverings. For
this reason it has been extended in \cite{Mu3},
in order to obtain the whole class of cyclic branched coverings of
2-bridge knots/links.

This new class of 4-coloured graphs $\Gt(n,p,q,c,c')$ depends on
five integer parameters, with $n,p>0$, $q\in \Z_{2p}$, $c,c'\in
\Z_n$, $\gcd(p,q)=1$ and $\gcd(n,c,c')=1$. Each $\Gt(n,p,q,c,c')$
is defined by the following four fixed-point-free involutions on
the set $V=\Z_n\times \Z_{ 2p}$: $$
\begin{array}{ll}
\tilde{\iota}_0=\iota_0,& \qquad
\tilde{\iota}_1(i,j)=\big(i+c'\eta (j),1-j\big),\\
\tilde{\iota}_2=\iota_2, & \qquad
\tilde{\iota}_3=\iota_3.
\end{array}
$$
The graph $\Gt(n,p,q,c,c')$ represents a connected, orientable
pseudomanifold $\St(n,p,q,c,c')$.

\begin{proposition} \cite{Mu3}
The following properties hold:
\begin{description}
\item[(i)] $\Gt(n,p,q,c,c')$ represents a 3-manifold if and
only if either $p$ is even or at least one of the
conditions (1) $c=0$, (2) $c'=0$, (3) $c=(-1)^qc'$ is satisfied.
\item[(ii)] $\St(n,p,q,c,1)\cong S(n,p,q,c)$.
\item[(iii)]
$\St(n,p,q,0,c')\cong\St(n,p,q,c,0)\cong\St(n,1,1,-c',c')\cong
\S^3$.
\item[(iv)] $\St(n,p,q,c,c')\cong\St(n,p,q,c',c)$.
\end{description}
\end{proposition}

We give the connection between generalized Lins-Mandel manifolds
and cyclic branched coverings of 2-bridge knots/links.

\begin{theorem} \cite{Mu3}
The $3$-manifold $\St(n,p,q,c,c')$, with $ c\ne 0\ne c'$, is
homeomorphic to $M_{n,c',-c}(p/q)$. Therefore, the class of
generalized Lins-Mandel manifolds $\tilde { S}(n,p,q,c,c')$, with
$c\ne 0\ne c'$ and $p>1$, is precisely the class of all cyclic
coverings of $\S^3$ branched over the 2-bridge knots/links
(with the exception of the trivial link with two components).
\end{theorem}

Another family of 4-coloured graphs representing cyclic branched
coverings of 2-bridge knots is described in \cite{KoSo}.

\section{Fundamental groups}

A presentation of the fundamental group of $M_{n,k}(\a, \b)$ can
be obtained from the coloured graph construction. Define the
following words in $x_1, \ldots , x_n$ with subscripts mod $n$: $$
Q_i = \prod_{j=0}^{n'-1} x_{i - j k}, \quad 1 \le i \le \gcd(n,k)
; \qquad \quad Q_i{'} = \prod_{j=0}^{\a-1}  x_{i+ s_j }^{e_j},
\quad 1\le i\le n; $$ with $n' = n / \gcd(n,k)$, $e_j=-\mu(2j\b)$
and $$ s_j  = \left\{
\begin{array}{ll}
\displaystyle
-k \sum_{h=1}^j \mu (\a + 2 \b - 2 h \b)-\sum_{h=1}^j\mu
(\a+\b-2h\b) &
\mbox{   if  } e_j = +1 \\
\displaystyle
-k \sum_{h=1}^j \mu (\a + 2 \b - 2 h \b)-\sum_{h=1}^j\mu
(\a+\b-2h\b) + k
&\mbox{  if  } e_j = -1
\end{array} \right.  ;
$$ where $\mu:\Z_{2\a}\to\{-1,+1\}$ is the map defined in
Section~6.

\begin{theorem} \cite{Mu3}
The fundamental group of $M_{n,k}(\a, \b)$ has the following
presentation $$ \pi_1 (M_{n,k}(\a, \b)) = \langle x_1, \ldots, x_n
\, |  \, Q_i = 1,  \,  1 \le i \le \gcd(n,k), \quad Q_i{'}= 1,
\,1\le i\le n \rangle ,$$ where $Q_i$ and $Q_i{'}$ are as above.
\end{theorem}

Since the cyclic branched coverings of 2-bridge knots have a
cyclic homeomorphism, it is natural to wonder about some cyclic
presentations of their fundamental groups. We recall that a finite
balanced presentation of a group $G\cong\langle
x_1,\ldots,x_n\vert r_1,\ldots,r_n \rangle$ is said to be a {\it
cyclic presentation\/} if there exists a word $w$ in the free
group $F_n$ generated by $x_1,\ldots,x_n$ such that the relators
of the presentation are $r_k=\theta_n^{k-1}(w)$, $k=1,\ldots,n$,
where $\theta_n :F_n\to F_n$ denotes the automorphism defined by
$\theta_n (x_i)=x_{i+1}$ (subscripts mod $n$), $i=1,\ldots,n$. Let
us denote this cyclic presentation (and the related group) by the
symbol $G_n(w)$, so that: $$ G_n(w) = \langle x_1, x_2, \ldots,
x_n \vert w, \theta_n(w), \ldots, \theta_n^{n-1}(w) \rangle. $$ A
group is said to be cyclically presented if it admits a cyclic
presentation. The polynomial associated with the cyclic
presentation $G_n (w)$ is given by $$ f_w (t)= \sum_{i=1}^{n}
a_i t^{i-1}, $$ where $a_i$ is the exponent sum of $x_i$ in $w$.
For the theory of cyclically presented groups we refer to \cite {Jo}.

Two different cyclic presentations for the fundamental groups of the
manifolds $M_n (\a/\b)$ are obtained in \cite{Mi} and
\cite{MuV}. Remark that explicit cyclic presentations different
from the above are listed in the Appendix of \cite{CRS}, for
2-bridge knots up to nine crossings.

{}From the polyhedral construction of $M_n (\a / \b)$ the following
presentation of the fundamental group holds. Denote by $$ R_{\a /
\b} (x_1, x_2, \ldots, x_n) = x_1 x_{1+s_1}^{-1} x_{1+s_2}
x_{1+s_3}^{-1} \cdots x_{1+s_{\a-1}}^{(-1)^{\a-1}}  , $$ with $$
s_j = s_j (\a, \b) = \sum_{i=1}^{j} (-1)^{\left[ i \b^{-1} / \a
\right]} , $$ where $\b^{-1}$ is the inverse of the element $\b$
in $\Z_{2\a}$ and $\left[ x \right]$ denotes the integral part of
$x$.

\begin{theorem} \label{MinkusPresentation}
\cite[Theorem~10]{Mi}
The fundamental group of $M_n(\a,\b)$
admits the presentation:
$$
 \langle x_1, \dots , x_n  |
R_{\a/\b} (x_i, \dots , x_{i+n-1}) = 1, \quad i=1, \dots ,
n \rangle \eqno(1)
$$
if $\a$ is odd,
$$
 \langle x_1, \dots , x_n , y | x_n=1,
R_{\a/\b} (x_i, \dots , x_{i+n-1}) = y, \, i=1, \dots , n \rangle
\eqno(2)
$$
if $\a$ is even.
\end{theorem}

Remark that (1) is a cyclic presentation and (2) is a non-cyclic one.

There is a nice relation between the Alexander polynomial
of a 2-bridge knot and the polynomial associated with
the above cyclic presentation.

\begin{proposition} \label{proposition D2} \cite[Theorem~11]{Mi}
The Alexander polynomial of the knot $\bt (\a, \b)$ is equal to
the polynomial associated with the cyclic presentation (1),
up to units of $\Z [t, t^{-1}]$.
\end{proposition}

According to \cite[Remark~4]{Mi}, this property holds for
a wider class of knots and cyclic presentations of their
cyclic branched coverings.

The following cyclic presentation for $\p_1(M_n(\a/\b))$, when
$\a$ is odd, arises from the surgery description of the manifold as
Takahashi manifold (see Section~5).

\begin{theorem} \label{cyclicpresentation} \cite{MuV}
Let $M_n(\a/\b)$ be the $n$-fold cyclic branched covering of the
2-bridge knot ${\bf b}(\a/\b)$, with $\a/\b=[-2q_1,2s_1, \ldots,
-2q_m, 2s_m]$. Then its fundamental group has the following cyclic
presentation:
$$
\pi_1 (M_n (\a/\b)) \, = \, \langle x_1, \dots ,
x_n  | w_{\a/\b} (x_i, \dots , x_{i+n-1}) = 1, \quad i=1, \dots ,
n \rangle ,
$$
where $$ w_{\a/\b} (x_i, \dots , x_{i+n-1}) =
b_{i+1,m}^{-s_m} d_{i+1,m} b_{i,m}^{s_m} ,$$
for $i=1,\ldots,n$ (subscripts mod $n$). The right part of these formulas are defined
by
the recurrent rule $$ d_{k,j} = b_{k,j-1}^{-s_{j-1}} d_{k,j-1}
b_{k-1,j-1}^{s_{j-1}}, \qquad b_{k,j} = d_{k,j}^{q_j} b_{k,j-1}
d_{k+1,j}^{-q_j}, \qquad j=2, \dots , m $$ and $$ b_{k,1} =
d_{k,1}^{q_1} d_{k+1,1}^{-q_1}, $$ where $x_k = d_{k,1}$, for
$k=1, \dots , n$.
\end{theorem}

We will illustrate the result obtained for $m=1$ and $m=2$.

\smallskip

If $m=1$, then $\a/\b = -2q + 1/(2s)$. Observe that this case
corresponds to Takahashi manifolds and was discussed in
\cite{K1KV, KV}. We get $$ \pi_1 (M_{n} ( -2q  + \frac{1}{2s} ) )
= \langle x_1, \dots , x_n \, | \, (x_k^q x_{k+1}^{-q})^{-s} x_k
(x_{k-1}^q x_k^{-q})^s = 1,
 \quad k=1,\dots ,n  \rangle .
$$

For example, if $q=-1$ and $s=1$ then $\a/\b=5/2$, and $\bt (5/2)$
is the figure-eight knot $4_1$ \cite{BZ}. So, its $n$-fold cyclic
branched covering has the fundamental group with the following
cyclic presentation $$ \pi_1 (M_n (5/2)) \, = \, \langle x_1,
\dots , x_n \, | \, \, x_{k+1}^{-1} x_k^2 x_{k-1}^{-1} x_k  = 1 ,
\quad k=1, \dots ,n \rangle $$ (compare with \cite{CRS, K1KV,
KV}).

If $m=2$ then $\bt (\a/\b)$ has Conway parameters $[-2q_1, 2s_1,
-2q_2, 2s_2]$ and $ \pi_1 (M_{n}(\a/\b))$ has the following
presentation: $$ \langle x_{1}, \dots , x_{n} \quad | \quad
w_{\a/\b} (x_{k-2}, x_{k-1}, x_{k}, x_{k+1}, x_{k+1}) = 1, \qquad
k=1, \dots , n \rangle , $$ where $$ w_{\a/\b} (x_{k-2}, x_{k-1},
x_{k}, x_{k+1}, x_{k+2}) = $$ $$ \left[ \left[ (x_k^{q_1}
x_{k+1}^{-q_1})^{-s_1} x_k (x_{k-1}^{q_1} x_{k}^{-q_1})^{s_1}
\right]^{q_2}
 x_k^{q_1} x_{k+1}^{-q_1}
\left[ (x_{k+1}^{q_1} x_{k+2}^{-q_1})^{-s_1} x_{k+1}
(x_{k}^{q_1}
x_{k+1}^{-q_1})^{s_1} \right]^{-q_2} \right]^{-s_2}
$$
$$
\cdot
(x_k^{q_1} x_{k+1}^{-q_1})^{-s_1} x_k (x_{k-1}^{q_1}
x_{k}^{-q_1})^{s_1}
$$
$$
\cdot \left[ \left[
(x_{k-1}^{q_1}
x_{k}^{-q_1})^{-s_1} x_{k-1} (x_{k-2}^{q_1}
x_{k-1}^{-q_1})^{s_1}
\right]^{q_2}
 x_{k-1}^{q_1} x_{k}^{-q_1}
\left[ (x_{k}^{q_1} x_{k+1}^{-q_1})^{-s_1} x_{k}
(x_{k-1}^{q_1}
x_{k}^{-q_1})^{s_1} \right]^{-q_2} \right]^{s_2} .
$$

For example, if $q_1=q_2=-1$ and $s_1=s_2=1$ then $\a/\b=29/12$,
that corresponds to the knot $8_{12}$. So, its $n$-fold cyclic
branched covering has the fundamental group with the following
cyclic presentation:
\begin{eqnarray}
\langle x_1, \dots, x_n \, | \, &
x_{k+1}^{-1} x_{k} x_{k+1}^{-2} x_{k+2} x_{k+1}^{-1} x_{k}
x_{k+1}^{-1} x_{k}^{2} x_{k-1}^{-1} x_{k} x_{k+1}^{-1}
x_{k}^{2} x_{k-1}^{-1}
&  \nonumber \\ &
\cdot
x_{k} x_{k-1}^{-1} x_{k-2} x_{k-1}^{-2} x_{k}
x_{k-1}^{-1} x_{k}
x_{k+1}^{-1} x_k^2 x_{k-1}^{-1} x_k = 1, &  k=1, \dots
,n \rangle .
\nonumber
\end{eqnarray}

Another cyclic presentation for the fundamental group of every
cyclic branched covering of a 2-bridge knots can be obtained
from the Heegaard diagram construction of \cite{Du}, since the
Dunwoody family contains these ma\-ni\-folds (see Theorem
\ref{Theorem
Dunwoody}). Moreover, these presentations arise from Heegaard
diagrams and are therefore geometric.

We point out that the problem of determining whether a balanced
presentation of a group is {\it geometric} -- i.e., induced by a
Heegaard
diagram of a closed orientable 3-manifold -- is of considerable
interest within geometric topology  \cite{Ne, OS1, OS2, OS3, St}.

\begin{corollary} \label{Corollary 6} \cite{GM}
The fundamental group of every branched cyclic
covering of a 2-bridge knot admits a cyclic presentation which is
geometric.
\end{corollary}

\section{Homology}

In this section we present the homology groups of some classes of
cyclic branched coverings of 2-bridge knots/links. Recall
(see, for example \cite[p.~71]{Ka}) that when $H_1 (M_n (\a/\b))$
is a finite abelian group, the order of this group is given by
the absolute value $| \Pi_{k=1}^{n} \Delta_{(\a,\b)} (\zeta^k) |$
where $\Delta_{(\a,\b)} (t)$ is the Alexander polynomial of $\bt
(\a,\b)$ and $\zeta$ is an $n$-th primitive root of unity.

As is well known, the manifold $M_n(\a/1)$ is homeomorphic to the
Brieskorn manifold $M(n,\a,2)$. Thus, in this case, the homology
groups can be obtained from \cite{Ra}.

\smallskip

If $\a$ is odd, then the $n$-fold cyclic branched covering is
unique and the homology groups are the following.

\begin{proposition} \cite{Ca4}
If $\a$ is odd then: $$H_1(M_n(\a/1))\cong \left\{
\begin{array}{ll}
\Z^{d-1}\oplus \Z_{\a/d}&\mbox{   if  }n\mbox{ is
even}\\\Z_2^{d-1}&\mbox{   if  } n\mbox{ is odd}\,
\end{array} \right.$$
where $d=\gcd(n,\a)$.
\end{proposition}

If $\a$ is even the situation is rather more complicated, since
$k$ can assume any value in $\Z_n-\{0\}$.

\begin{proposition} \cite{Mu3}
Let $\a$ be even. Then: $$H_1(M_{n,k}(\a/1))\cong \left\{
\begin{array}{ll}
\Z^{d-m}\oplus \Z_{a}^m & \mbox{  if }h=1 \\ \Z^{d-h+1-m}\oplus
\Z_{a}^{m-h+1}\oplus \Z_{ab}^{h-2}\oplus \Z_{hab}&\mbox{  if
}1<h<m+1\\ \Z^{d-h+1-m}\oplus \Z_{b}^{h-1-m}\oplus
\Z_{ab}^{m-1}\oplus \Z_{ hab}&\text{  if }h\ge m+1
\end{array}
\right. .$$ where $s=\gcd(n,k)$, $d=\gcd(n,\a(k+1)/2)$,
$h=\gcd(n,k+1)$, $m=\gcd(d,s)$, $a=nm/(sd)$ and $b=\a h/(2d)$.
\end{proposition}

As a consequence we have the homology groups in the interesting cases
$\gcd(n,k)=1$ (i.e., meridian-cyclic coverings).

\begin{corollary}
Let $\a$ be even and $\gcd(n,k)=1$. Then:
$$H_1(M_{n,k}(\a/1))\cong \left\{
\begin{array}{ll}
\Z^{d-1} \oplus \Z_{a} & \mbox{  if  }h=1\ \\ \Z^{d-h}\oplus
\Z^{h-2}_{b}\oplus \Z_{hab}& \mbox{ if }h>1
\end{array} ,\right.$$
where $d=\gcd(n,\a(k+1)/2)$, $h=\gcd(n,k+1)$, $a=n/d$ and $b=\a
h/(2d)$.
\end{corollary}

The previous corollary contains, as a particular case, the
homology groups of the Brieskorn manifolds $M(n,\a,2)\cong
M_n(\a/1)$ (see also \cite{Ca4}).

\begin{corollary} If $\a$ is even, then: $$H_1(M_n(\a/1))\cong
\left\{
\begin{array}{ll}
\Z^{d-1}\oplus \Z_{ n/d} & \mbox{   if
}n \mbox{ is odd}
\\ \Z^{d-2}\oplus \Z_{2n\a/d^2} & \mbox{
if  } n
\mbox{ is even}
\end{array},
\right. $$ where $d=\gcd(n,\a)$.
\end{corollary}

The homology of cyclic branched coverings of 2-bridge knots of
genus one has been obtained in \cite{Mu3} and \cite{CRS}, using an
algorithm of Fox \cite{Fo2}. Recall that the 2-bridge knot
$\bt(\a,\b)$ has genus one if and only if $\b/2$ divides
$(\a-1)/4$ when $\a\cong 1$ mod $4$ and $\b/2$ divides $(\a+1)/4$
when $\a\cong 3$ mod $4$ (up to equivalence we can assume that $\b$
is even for any 2-bridge knot $\bt(\a,\b)$).

\begin{proposition} \cite{Mu3,CRS}
Let $\bt(\a,\b)$ be a 2-bridge knot of genus one. Then:
$$
H_1(M_n(\a/\b))\cong\left\{
\begin{array}{ll}
\Z_{\a\vert A'(n)\vert}\oplus \Z_{\vert A'(n)\vert}&
\mbox{  if  } n\mbox{ is even}\\
\Z_{\vert A''(n)\vert}\oplus \Z_{ \vert A^{\prime\prime}(n)\vert}&
\mbox{  if  }n\mbox{ is odd}
\end{array}
\right. ,
$$
where $A'(n),A''(n)$ are the integers defined by:
$$
A'(1)=1,\, A'(2)=1, \,A'(n+2)=A'(n+1)- h A'(n),
$$
$$A''(1)=1,\,A'{}'(2)=1-2h ,\,A'{}'
(n+2)=A''(n+1)-h A''(n),
$$
and
$$
h= \left\{
\begin{array}{ll}
(1-\a)/4 &\mbox{  if  } \a\cong 1\mbox{ mod }4
\\
(1+\a)/4 &\mbox{  if  } \a\cong 3\mbox{ mod }4
\end{array}
\right. .
$$
Thus, the homology does not depend on $\b$.
\end{proposition}

In some particular cases, explicit formulae are obtained
by J.~Minkus.

\begin{proposition} \cite[Corollary 11.2]{Mi} \begin{description}
\item[(i)] If $n>1$ is even then
$H_1(M_n((2n\b\pm 1)/\b))\cong\Z_{2n\b\pm 1}.$
\item[(ii)] If $n>0$ is odd then $M_n((2n\b\pm 1)/\b)$ is
a homology sphere.
\end{description}
\end{proposition}

Cyclic branched coverings of the Whitehead link $\bt(8,3)$ have
been intensively studied (see \cite{CP,HKM1,Zi1}).

\begin{proposition} \label{Whitehead-monodromy} \cite{HKM1}
Let $n\ge 3$ and let $\gcd(n,k)=1$, then:
$$H_1(M_{n,k}(8/3))\cong
\left\{
\begin{array}{ll}
\Z_{n/6}\oplus \Z_{n/2}\oplus \Z_{12n}& \mbox{  if
}n\equiv 0 \mbox{ mod }6\\
\Z_{n/2}\oplus \Z_{n/2}\oplus \Z_{4n}& \mbox{  if
}n\equiv \pm 2 \mbox{ mod }6\\
\Z_{n/3}\oplus \Z_{n}\oplus \Z_{3n}& \mbox{  if
}n\equiv 3 \mbox{ mod }6\\
\Z_{n}\oplus \Z_{n}\oplus \Z_{n}& \mbox{  if
}n\equiv \pm 1 \mbox{ mod }6
\end{array}
\right. .$$
\end{proposition}

\smallskip


For $n, \a \le 9$ and $gcd(n,k)=1$,
the homology groups of $M_{n,k}(\a/\b) \cong S(n,\a,\b,-k)$
are listed in \cite{Li} and \cite{LM}.

\section{Decomposition of singly-cyclic coverings}

In this section we prove that each singly-cyclic branched covering
of a 2-bridge link is the composition of a meridian-cyclic
branched covering of a certain link ${\cal L} (d, \a / \b)$
described below and a cyclic branched covering of a trivial knot.
This gives a generalization of a result obtained in \cite{CP} for
the case of the Whitehead link.

We can always assume that
$\a / \b  = \left[ 2a_1, - 2b_1, \dots, 2a_l, - 2b_l\right]$
if $\bt (\a, \b)$ is a knot, and
$\a / \b = \left[ 2a_1, -2b_1, \dots, 2a_l \right]$
if $\bt (\a, \b)$ is a link.
For any $\a / \b$
and $d \ge 1$, define a link ${\cal L} (d, \a /
\b)\subset\S^3$ as in Figure~13, where  $a_i$ and $b_i$ denote
numbers
of half-twists (in the vertical direction) in the corresponding
boxes,
and the fragment ``in degree'' $d$ must be repeated $d$ times.

\begin{figure}
 \begin{center}
 \includegraphics*[totalheight=14cm]{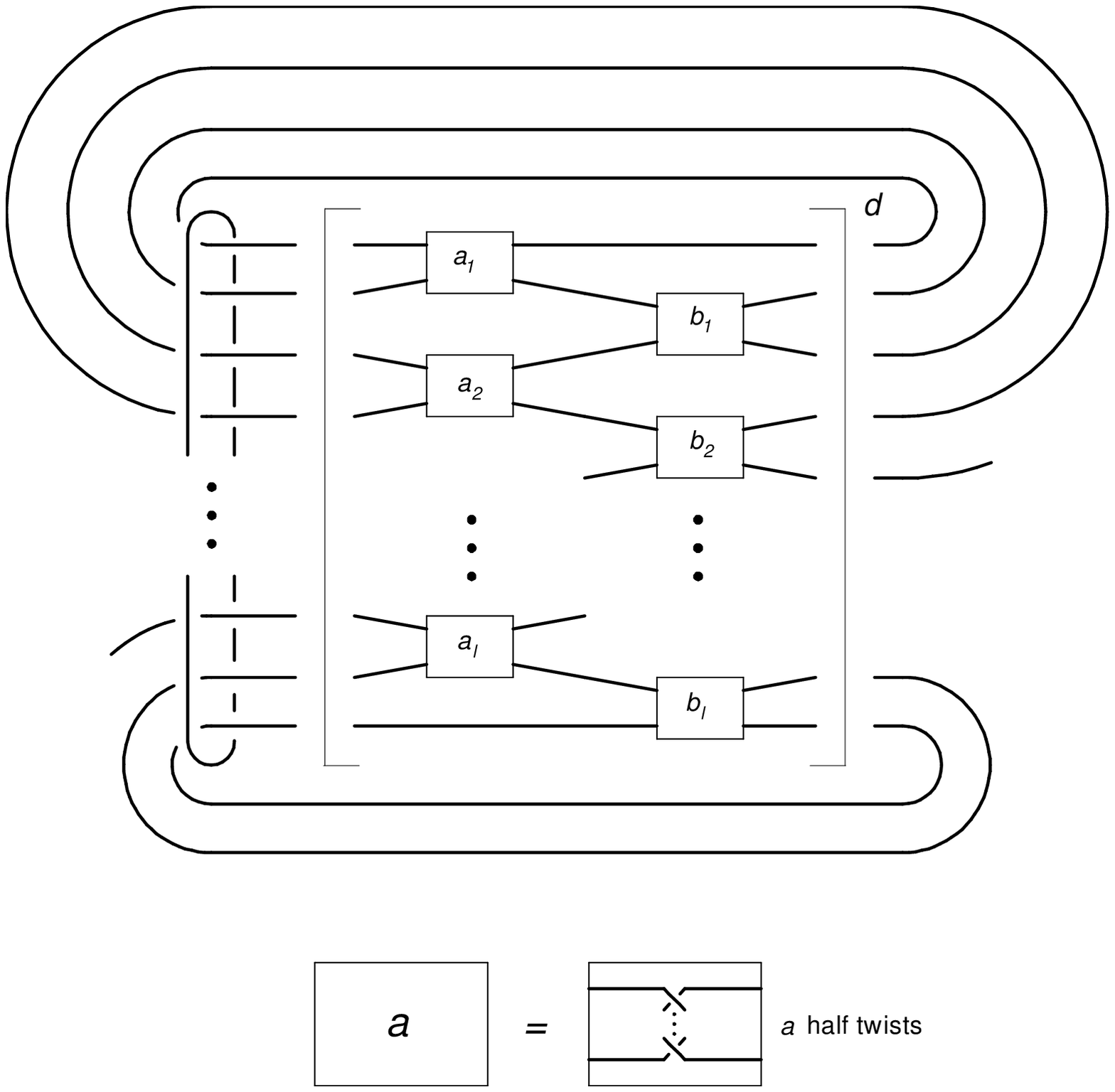}
 \end{center}
 \caption{ The link ${\cal L}(d, \a / \b)$. }
 \label{Fig. 14}
\bigskip\bigskip
\end{figure}

The link ${\cal L} (d, \a / \b)$ has $1+\gcd(d,l)$ components,
where $$l=\mbox
{lk}(\bt(\a,\b))=\sum_{h=1}^{\a/2}(-1)^{\left[\frac{(2h-1)\b}{\a}\right]}$$
is the linking number of $\bt(\a,\b)$ \cite[Remark 12.7]{BZ}
($\left[ x
\right]$ denotes the integral part of $x$). For example $\mbox
{lk}(\bt(8,3))=0$ and therefore ${\cal L} (d,8/3)$ has $d+1$
components.

Observe that, roughly speaking, each 2-bridge link can be obtained
as a quotient of a suitable 2-bridge knot/link by an involution
whose axis does not intersect the knot/link (for example, see the
2-periodic presentations of 2-bridge knots/links in \cite{BZ}).
Hence, any 2-bridge link $\bt (\a, \b)$ can be presented in the
form ${\cal L} (1, \a_1 / \b)$, where $\a_1 = \a /2$. For example,
if $\bt (8/3)$ is the Whitehead link, then $\bt (4/3)$ is
equivalent to $\bt (4/1)$, and $4/1 = [2 a_1]$ with $a_1=2$.
Hence we get the presentation of the Whitehead link in the form
${\cal L} (1, 4/1)$ (see, for example, Figure~3.4 in \cite{CP}).

Denote by ${\cal L}_{m} (d, \a / \b)$ the orbifold with underlying
space $\S^3$ and singular set ${\cal L} (d, \a / \b)$, with all
branched indices equal to $m$. Moreover, denote by $\bt_{n,n'}(\a,
\b)$ the orbifold having $\S^3$ as underlying space and the
2-bridge link $\bt( \a , \b)$ as singular set, with branched
indices $n$ and $n'$ respectively for the two components.

\begin{theorem} \label{singly}
Let $M_{n,k}(\a / \b)$ be a singly-cyclic covering of the
2-bridge link $\bt (\a, \b)$, and $\a_1 = \a /2$.
Suppose $d = \gcd(n,k)$. Then the following diagram
\smallskip
\begin{center}
  \unitlength=0.5mm
  \begin{picture}(90,75)(0,10)
  \put(40,80){\makebox(0,0)[cc]{$M_{n,k}(\a /\b)$}}
  \put(50,70){\vector(2,-1){30}}
  \put(40,70){\vector(0,-1){50}}
  \put(68,68){\makebox(0,0)[cc]{$n/d$}}
  \put(35,45){\makebox(0,0)[cc]{$n$}}
  \put(80,45){\makebox(0,0)[cc]{${\cal L}_{n/d}(d, \a_1 / \b)$}}
  \put(80,35){\vector(-2,-1){30}}
  \put(68,25){\makebox(0,0)[cc]{$d$}}
  \put(40,10){\makebox(0,0)[cc]{$\bt_{n,n/d}(\a, \b)$}}
  \end{picture}
\end{center}
\smallskip
is commutative. Moreover, the $n/d$--covering is a meridian-cyclic
covering of the link ${\cal L}(d, \a_1 / \b)$ and the $d$-covering
is the $d$-fold cyclic branched covering of a component of $\bt
(\a, \b)$ (which is a trivial knot).
\end{theorem}

\begin{proof}
The fundamental group $\Gamma$ of the orbifold $\bt_{n, n/d} (\a /
\b)$ admits the presentation $$ \Gamma = \langle \mu_1, \mu_2 \ |
\ \mu_1^n = \mu_2^{n/d} = 1, \, w(\mu_1, \mu_2)=1 \rangle $$ where
$\mu_1$ and $\mu_2$ are the homotopy classes of two meridians
$m_1$ and $m_2$ around the components $K_1$ and $K_2$ of the link,
and $w(\mu_1, \mu_2)$ is the relation deriving from the standard
presentation of the group of $\bt (\a, \b)$. From the definition
of a singly-cyclic covering we have $\pi_1 (M_{n,k}) = \ker
(\varphi)$ for the epimorphism $\varphi : \Gamma \to \Z_n =
\langle \gamma | \gamma^n = 1 \rangle$ defined by $\varphi (\mu_1)
= \gamma$ and $\varphi (\mu_2) = \gamma^k$. Consider the subgroup
$\Z_d \triangleleft \Z_n$ such that $\Z_d = \langle \delta \, | \,
\delta^d=1\rangle$, where $\delta = \gamma^{n/d}$. Let $\theta :
\Gamma \to \Z_d$ be the epimorphism defined by $\theta(\mu_1) =
\delta$ and $\theta(\mu_2)=1$. This epimorphism induces the
$d$-fold cyclic covering $\Theta : {\cal O} \to \bt_{n,n/k}
(\alpha, \beta)$ such that the axis of  the cyclic group action is
the component $K_1$ of $\bt (\a , \b)$ corresponding to the
meridian $\mu_1$. Therefore, the underlying space of the orbifold
$\cal O$ is $\S^3$ and the singular set is $\Theta^{-1} (K_1) \cup
\Theta^{-1} (K_2)$. Obviously, $\Theta^{-1} (K_1)$ is a trivial
knot in $\S^3$ with singularity index $n/d$, and $\Theta^{-1}
(K_2)$ is a $d$-periodic knot/link which also has singularity
index $n/d$. Since $\bt (\a, \b)$ is equivalent to ${\cal L} (1,
\a_1 / \b)$ with $\a_1 = \a / 2$, then we have $\Theta^{-1}
(\bt (\a, \b)) = {\cal L} (d, \a_1 / \b)$. Using $\pi_1 ({\cal O})
= \ker (\theta) = \varphi^{-1} (\Z_d)$ and $\pi_1 (M_{n,k}) =
\ker(\varphi)$, the diagram of coverings is commutative. This
completes our proof.
\end{proof}






\bigskip\bigskip

\vspace{15 pt} {MICHELE MULAZZANI, Department of
Mathematics,
University of Bologna, I-40127 Bologna, ITALY, and
C.I.R.A.M.,
Bologna, ITALY. E-mail: mulazza@dm.unibo.it}

\vspace{15 pt} {ANDREI VESNIN, Sobolev Institute
of Mathematics, Novosibirsk, 630090 Russia. E-mail:
vesnin@math.nsc.ru}



\begin{thebibliography}{5}

\bibitem {Al}
Alexander, J.W.:
Note on Riemann spaces.
Bull. Amer. Math. Soc. {\bf 26} (1920), 370--373.

\bibitem {BKM}
Bandieri, P., Kim, A.C., Mulazzani, M.,:
On the cyclic coverings of the knot $5_2$.
Proc. Edinb. Math. Soc. {\bf 42} (1999), 575--587.

\bibitem {BH}
Birman, J.S., Hilden, H.M.:
Heegaard splittings of $\S^3$.
Trans. Amer. Math. Soc. {\bf 213} (1975), 315--352.


\bibitem {BZ}
Burde, G., Zieschang, H.:
Knots.
de Gruyter Studies in Mathematics, {\bf 5},
Berlin--New York, 1985.

\bibitem{CG1}
Casali, M.R., Grasselli, L.: Characterizing crystallizations among
Lins-Mandel $4$-coloured graphs. Rend.  Circolo Mat.  Palermo,
Serie II {\bf 18} (1988), 221--228.

\bibitem{CG2}
Casali, M.R., Grasselli, L.:
$2$-Symmetric crystallizations and $2$-fold branched
coverings of $\S^3$.
Discrete Math. {\bf 87} (1991), 9--22.

\bibitem{Ca1}
Cavicchioli, A.:
Lins-Mandel crystallizations.
Discrete Math. {\bf 57} (1985), 17--37.

\bibitem{Ca2}
Cavicchioli, A.:
A countable class of non-homeomorphic homology
spheres with Heegaard
genus two.
Geom. Dedicata {\bf 20} (1986), 345--348.

\bibitem{Ca3}
Cavicchioli, A.: Lins-Mandel $3$-manifolds and their groups:  a
simple proof of the homology sphere conjecture. Rend.  Circolo
Mat. Palermo, Serie II {\bf 18} (1988), 229--237.

\bibitem{Ca4}
Cavicchioli, A.:
On some properties of the groups $G(n,l)$.
Ann. Mat. Pura Appl. {\bf 151} (1988), 303--316.

\bibitem {CHK1}
Cavicchioli, A., Hegenbarth F., Kim, A.C.: A
geometric study of
Sieradski groups. Algebra Colloq. {\bf 5} (1998),
203--217.

\bibitem {CHK2}
Cavicchioli, A., Hegenbarth F., Kim, A.C.: On cyclic
branched
coverings of torus knots. J. of Geometry {\bf 64}
(1999), 55--66.

\bibitem {CHR}
Cavicchioli, A., Hegenbarth, F., Repovs, D.: On manifold spines
and cyclic presentations of groups. In: Knot theory. Banach Center
Publ. {\bf 42} (1998), 49--56.

\bibitem{CP}
Cavicchioli, A., Paoluzzi, L.: On certain classes of
hyperbolic
3-manifolds. Manuscripta Math. {\bf 101} (2000),
457--494.

\bibitem{CRS}
Cavicchioli, A., Ruini, B., Spaggiari, F.:
Cyclic branched coverings of 2-bridge knots.
Revista Mat. Compl. {\bf 12} (1999), 383--416.

\bibitem{CS}
Cavicchioli, A., Spaggiari, F.: The classification of
$3$-manifolds with spines related to Fibonacci groups. In: Lect.
Notes in Math, Springer Verlag {\bf 1509} (1992), 50--78.

\bibitem {Co}
Conway, J.H.: An enumeration of knots and links and some of their
related properties. In: Computational problems in Abstract
Algebra, Proc. Conf. Oxford 1967 (ed. J.Leech), Pergamon Press
(1969), 329--358.

\bibitem{Do}
Donati, A.: Lins-Mandel manifolds as branched
coverings of $
\S^3$. Discrete Math. {\bf 62} (1986), 21--27.

\bibitem{DG}
Donati, A., Grasselli, L.: Gruppo dei colori e
cristallizazioni
``normali'' degli spazi lenticolari. Boll. Un. Mat.
Ital. {\bf
16,1-A} (1982), 359--366.

\bibitem {Dn}
Dunbar, W.D.: Geometric Orbifolds.
Rev. Mat. Univ. Complutense Madr. {\bf 1} (1988), 67--99.

\bibitem {Du}
Dunwoody, M.J.: Cyclic presentations and
3-manifolds. In: Proc.
Inter. Conf., Groups-Korea '94 (A.C. Kim, D.L. Johnson eds.), Walter
de Gruyter,
Berlin-New York
(1995), 47--55.


\bibitem{FGG}
Ferri, M., Gagliardi, C., Grasselli, L.: A
graph-theoretical
representation of PL-manifolds -- A survey on
crystallizations. Aequationes
Math. {\bf 31} (1986), 121--141.

\bibitem{FM}
Fomenko, A.T., Matveev, S.V.: Constant energy
surfaces of
Hamiltonian systems, enumeration of
three-dimensional manifolds in
increasing order of complexity, and computation of
volumes of
closed hyperbolic manifolds. Russian Math. Surveys
(1) {\bf 43}
(1988), 3--24.

\bibitem{Fo1}
Fox, R.H.: Covering spaces with singularities. In:
Algebraic
Geometry and Topology, a simposium in honour of
S.Lefschetz,
Princeton Math. Series {\bf 12} (1957), 243--257.

\bibitem{Fo2}
Fox, R.H.: The homology characters of the cyclic
coverings of the
knots of genus one. Ann. of Math. {\bf 71} (1960),
187--196.

\bibitem{Fu}
Funcke, K.: Geschlecht von Knoten mit zwei Br\"ucken
und die
Faserbarkeit ihrer Aussenr\"aume. Math. Z. {\bf 159}
(1978), 3--24.

\bibitem{Gr}
Grasselli, L.: The groups $G(n,l)$ as fundamental
groups of Seifert fibered
homology spheres. Proc. ``Groups 1993 -
Galway/St.Andrews'' {\bf 1},
Cambridge Univ. Press (1994), 244--248.

\bibitem {GM}
Grasselli, L., Mulazzani, M.:
Genus one 1-bridge knots and Dunwoody manifolds.
Forum Math. {\bf 13} (2001), 379--397.

\bibitem {HKM1}
Helling, H., Kim, A.C., Mennicke, J.L.: Some
honey-combs in
hyperbolic 3-space. Commun. Algebra {\bf 23} (1995),
5169--5206.

\bibitem {HKM2}
Helling, H., Kim, A.C., Mennicke, J.L.: A geometric
study of
Fibonacci groups. J. Lie Theory {\bf 8} (1998), 1-23.

\bibitem{He}
Hempel, J.: 3-manifolds. Annals of Math. Studies,
vol. 86,
Princeton University Press, Princeton, N.J., 1976.

\bibitem {HLM1}
Hilden, H.M., Lozano, M.T., Montesinos-Amilibia,
J.M.:
The arithmeticity of the figure eight knot orbifolds.
In: `Topology 90' (B. Apanasov, W. D. Neumann, A. W.
Reid,
and L. Siebenmann, Eds.) Ohio State Univ., Math.
Research Inst.
Publ., Walter de Gruyter Ed., Berlin {\bf 1} (1992),
169--183.

\bibitem {HLM2}
Hilden, H.M., Lozano, M.T., Montesinos-Amilibia,
J.M.:
On the arithmetic $2$-bridge knots and link orbifolds
and a new knot invariant.
J. Knot Theory Ramifications {\bf 4} (1995), 81--114.

\bibitem {HLM3}
Hilden, H.M., Lozano, M.T., Montesinos-Amilibia,
J.M.:
Volumes and Chern--Simons invariants of cyclic
coverings over rational knots.
In: Proc. of The 37-th Tanigughi Symposium on
Topology and Teichm\"uller
spaces held in Finland, July 1995, ed. by S.~Kojima
et al..
World Scientific Publishing Co. (1996), 31--55.

\bibitem {HW}
Hilton, P.J., Wylie, S.: An introduction to algebraic topology -
Homology theory. Cambridge, 1960.

\bibitem{Ho}
Hodgson, C.:
Degeneration and regeneration of geometric
structures on three-manifold,
Ph. D. Thesis, Princeton University, 1986.

\bibitem {HoW}
Hodgson, C., Weeks, J.:
Symmetries, isometries and length spectra of closed
hyperbolic
three-manifolds.
Exper. Math. {\bf 3} (1994), 101--113.


\bibitem {Jo}
Johnson, D.L.: Topics in the theory of group presentations. London
Math. Soc. Lect. Note Ser., vol. 42, Cambridge Univ. Press,
Cambridge, U.K., 1980.

\bibitem{JT}
Johnson, D.L., Thomas, R.M.: The Cavicchioli groups
are pairwise
non-isomorphic. London Math. Soc. Lecture Note
Series {\bf 121} (1986),
220--222.


\bibitem {Ka}
Kawauchi, A.: A survey of knot theory. Birkhauser
Verlag, Basel,
1996.

\bibitem {Ki}
Kim, A.C.: On the Fibonacci group and related
topics. Contemp.
Math. {\bf 184} (1995), 231--235.

\bibitem {K1KV}
Kim, A.C., Kim, Y., Vesnin, A.: On a class of
cyclically presented
groups. In: Proc. Inter. Conf., Groups-Korea '98 (Y.G. Baik, D.L.
Johnson, A.C. Kim eds.),
Walter de
Gruyter, Berlin-New York (2000), 211--220.

\bibitem {KV}
Kim, A.C., Vesnin, A.: The fractional Fibonacci
groups and
manifolds. Sib. Math. J. {\bf 39} (1998), 655--664.

\bibitem {K2KV}
Kim, G., Kim, Y., Vesnin, A.: The knot $5_2$ and
cyclically
presented groups. J. Korean Math. Soc. {\bf 35}
(1998), 961--980.

\bibitem {KiY}
Kim, Y.: About some infinite family of 2-bridge knots and
3-manifolds. Internat. J. Math. \& Math. Sci. {\bf 24} (2000),
95--108.

\bibitem {KS} Kirby, R.C.; Scharlemann, M.G.:
Eight faces of the Poincar\'e homology 3-sphere.
Geometric
topology, Proc. Conf., Athens/Ga. (1979) 113--146.

\bibitem {KoSo} Ko, K.H., Song, H.J.:
On the crystallization of 3-manifolds associated with polyhedral
schemata.
J. Korean Math. Soc. {\bf 36} (1999), 447--458.

\bibitem{Li}
Lins, S.: Gems, computers and attractors for
3-manifolds. World
Scientific, 1995.

\bibitem{LM}
Lins, S., Mandel, A.: Graph-encoded $3$-manifolds.
Discrete Math.
{\bf 57} (1985), 261--284.

\bibitem{LMu}
Lins, S., Mulazzani, M.: Isomorphisms and
homeomorphisms of a class of graphs and spaces.
To appear in Aequationes
Math.

\bibitem {MR}
Maclachlan, C., Reid, A.W.: Generalised Fibonacci
manifolds.
Transform. Groups {\bf 2} (1997), 165--182.

\bibitem{MM}
Mayberry, J., Murasugi, K.: Torsion groups of
abelian coverings
of links. Trans. Amer. Math. Soc. {\bf 271} (1982),
143--173.

\bibitem {MeV0}
Mednykh, A., Vesnin, A.:
Hyperbolic volumes of Fibonacci manifolds.
Sib. Math. J. {\bf 36} (1995), 235--245.

\bibitem {MeV1}
Mednykh, A., Vesnin, A.:
Fibonacci manifolds as two-fold coverings over the
three-dimensional
sphere and the Meyerhoff~--~Neumann conjecture.
Sib. Math. J. {\bf 37} (1996), 461--467.

\bibitem {MeV3}
Mednykh, A., Vesnin, A.:
The Heegaard genus of three-dimensional hyperbolic manifolds
of small volume.
Sib. Math. J. {\bf 37} (1996), 893--897.

\bibitem {MeV2}
Mednykh, A., Vesnin, A.:
Visualization of the isometry group action on the
Fomenko--Matveev--Weeks
manifold.
J. of Lie Theory  {\bf 8} (1998), 51--66.

\bibitem{M}
Milnor, J.: On the $3$-dimensional Brieskorn
manifolds $M(p,r,q)$.
In: Knots, groups and 3-manifolds (ed. by L. P.
Neuwirth) Ann.
Math. Studies {\bf 84} (1975), 175--225.

\bibitem {Mi}
Minkus, J.: The branched cyclic coverings of 2 bridge knots and
links. Mem. Amer. Math. Soc. {\bf 35} Nr. 255 (1982), 1--68.

\bibitem {Mo2}
Montesinos, J.M.: Representing 3-manifolds by a
universal
branching set. Math. Proc. Camb. Philos. Soc. {\bf
94} (1983),
109--123.

\bibitem{MB}
Morgan, J.W., Bass, H.: The Smith conjecture.
Academic
Press, Inc., 1984.

\bibitem {Mu2}
Mulazzani, M.:
Lins-Mandel graphs representing 3-manifolds.
Discrete Math. {\bf 140} (1995), 107--118.

\bibitem {Mu3}
Mulazzani, M.:
All Lins-Mandel spaces are branched cyclic coverings
of $\S^3$.
J. Knot Theory Ramifications {\bf 5} (1996),
239--263.

\bibitem{Mu1}
Mulazzani, M.:
A "universal" class of 4-coloured graphs.
Rev. Mat. Univ. Complutense Madr. {\bf 9} (1996), 165--195.

\bibitem {Mu4}
Mulazzani, M.: On $p$-symmetric Heegaard splittings.
J. Knot Theory Ramifications {\bf 9} (2000), 1059--1067.


\bibitem {MuV}
Mulazzani, M., Vesnin, A.:
Generalized Takahashi manifolds.
xxx archive preprint math.GT/0106137.

\bibitem{NO}
Negami, S., Okita, K.: The splittability and triviality of 3-bridge
links.
Trans. Amer. Math. Soc. {\bf 289}
(1985), 253--280.

\bibitem {Ne} Neuwirth, L.: An algorithm for the construction of
3-manifolds from 2-complexes. Proc. Camb. Philos. Soc. {\bf 64}
(1968), 603--613.

\bibitem {OS1}
Osborne, R.P., Stevens, R.S.: Group presentations corresponding to
spines of 3-manifolds I. Amer. J. Math. {\bf 96} (1974), 454--471.

\bibitem {OS2}
Osborne, R.P., Stevens, R.S.: Group presentations corresponding to
spines of 3-manifolds II. Trans. Amer. Math. Soc. {\bf 234}
(1977), 213--243.

\bibitem {OS3}
Osborne, R.P., Stevens, R.S.: Group presentations corresponding to
spines of 3-manifolds III. Trans. Amer. Math. Soc. {\bf 234}
(1977), 245--251.


\bibitem{Pe}
Pezzana, M.: Diagrammi di Heegaard e triangolazione contratta,
Boll. Un. Mat. Ital. {\bf 12}(1975), 98--105.

\bibitem{Ra}
Randell, R.C.: The homology of generalized Brieskorn manifolds.
Topology {\bf 14} (1975), 347--355.

\bibitem {RV}
Rasskazov, A., Vesnin, A.: Isometries of the
Hyperbolic
Fibonacci Manifolds. Siberian Math. J. {\bf 40}
(1999), 9--22.

\bibitem {Ro}
Rolfsen, D.: Knots and Links. Publish or Perish
Inc., Berkeley
Ca., 1976.

\bibitem {RS}
Ruini, B., Spaggiari, F.: On the structure of
Takahashi manifolds.
Tsukuba J. Math. {\bf 22} (1998), 723--739.
Corrigendum to ``On the structure of Takahashi manifolds''.
Tsukuba J. Math. {\bf 24} (2000), 433--434.

\bibitem {Sa}
Sakuma, M.: The geometries of spherical Montesinos
links. Kobe J.
Math. {\bf 7} (1990), 167--190.

\bibitem{Sc}
Schubert, H.: Knoten mit zwei Br\"ucken. Math. Z.
{\bf 65}
(1956), 133--170.

\bibitem{ST}
Seifert, H., Threlfall, W.:
A textbook of topology.
Academic Press (English reprint), 1980.

\bibitem {Si}
Sieradski, A.J.:
Combinatorial squashings, 3-manifolds, and the third
homotopy of groups.
Invent. Math. {\bf 84} (1986), 121--139.


\bibitem {SH}
Song, H.J., Hong, W.C.:
Extended Heegaard diagrams of cyclic branched coverings of $\S^3$
over 2-bridge links. Preprint.


\bibitem {SK}
Song, H.J., Kim, S.H.: Dunwoody 3-manifolds and
$(1,1)$-decomposible knots. In: Proceedings of Workshop in Pure
Mathematics, (Jongsu Kim and Sungbok Hom Ed.), ``Geometry and
Topology'', Vol. 19 (2000), 193--211.

\bibitem {St}
Stevens, R.S.: Classification of 3-manifolds with certain spines.
Trans. Amer. Math. Soc. {\bf 205} (1975), 151--166.

\bibitem {Ta}
Takahashi, M.:
On the presentations of the fundamental groups of
3-manifolds.
Tsukuba J. Math. {\bf 13} (1989), 175-189.

\bibitem {Th}
Thurston, W.: The geometry and topology of 3-manifolds. Notes
1976-1978, Princeton University Press.


\bibitem {Zi1}
Zimmermann, B.:
On the Hantzsche--Wendt manifold.
Monatsh. Math. {\bf 110} (1990), 321--327.

\bibitem{Zi5}
Zimmermann, B.:
Finite group actions on handlebodies and equivariant Heegaard genus
for 3-manifolds.
Topology and its Appl. {\bf 43} (1992), 263--274.

\bibitem {Zi2}
Zimmermann, B.:
On cyclic branched coverings of hyperbolic links.
Topology  Appl. {\bf 65} (1995), 287--294.

\bibitem {Zi3}
Zimmermann, B.:
Genus action of finite groups on 3-manifolds.
Michigan Math. J. {\bf 43} (1996), 593--610.

\bibitem {Zi4}
Zimmermann, B.: Determining knots and links by cyclic branched
coverings. Geom. Dedicata {\bf 66} (1997), 149--157.


\end{thebibliography}
\end{document}